\documentclass[a4paper,12pt,reqno]{amsart}
\usepackage{graphics}
\newtheorem{lem}{Lemma}

\newtheorem{theo}{Theorem}

\newtheorem{cor}{Corollary}

\numberwithin{equation}{section}

\newcommand{\sgn}{\operatorname{sgn}}

\newcommand{\id}{\operatorname{id}}

\begin{document}

\title[Halved monotone triangles]{An operator formula for the number of halved monotone triangles with prescribed bottom row}

\author[Ilse Fischer]{\box\Adr}

\newbox\Adr
\setbox\Adr\vbox{ \centerline{ \large Ilse Fischer} \vspace{0.3cm}
\centerline{Fakult\"at f\"ur Mathematik, Universit\"at Wien}
\centerline{Nordbergstrasse 15, A-1090 Wien, Austria}
\centerline{E-mail: {\tt Ilse.Fischer@univie.ac.at}} }

\begin{abstract}
Monotone triangles are certain triangular arrays of integers, which correspond to 
$n \times n$ alternating sign matrices when prescribing $(1,2,\ldots,n)$ 
as bottom row of the monotone triangle. In this article we define
halved monotone triangles, a specialization of which correspond to
vertically symmetric alternating sign matrices. We derive an
operator formula for the number of halved monotone triangles with
prescribed bottom row which is analogous to our operator formula
for the number of ordinary monotone triangles~\cite{fischer}.
\end{abstract}

\maketitle

\section{Introduction}

Alternating sign matrices and equivalent objects such as fully
packed loop configurations, the six vertex model and monotone
triangles are nowadays a rich source for intriguing problems on
which combinatorialists can test their various enumeration
methods. This article is another contribution in this respect.

In \cite{fischer} we gave a formula for the number of monotone
triangles with prescribed bottom row. Strikingly this formula
involves shift operators which are applied to a simple
multivariate polynomial. It is an example of a new type of
enumeration formula combinatorialists can possibly make use of when
answering their enumeration problems. Subsequently, our formula
enabled us to give a new proof of the refined alternating sign
matrix theorem~\cite{fischer1}, which was first proved by
Zeilberger~\cite{zeilberger2}. Here, we present a second example
of such an operator formula. This new formula gives the number of
halved monotone triangles with prescribed bottom row, a notion to be defined
below.

To keep the treatment self-contained, we recall the basic
definitions. An alternating sign matrix is a square matrix with
$0$s, $1$s and $-1$s as entries such that the row- and columnsums
are $1$ and the non-zero entries of each row and of each column
alternate in sign. Thus,
$$
\left(
\begin{array}{rrrrrrr}
0 & 0 & 0 & 1 & 0 & 0 & 0 \\
0 & 1 & 0 & -1 & 0& 1  & 0 \\
1 & -1 & 0 & 1 & 0 & -1 & 1 \\
0 & 0 &  1 & -1 & 1 & 0 & 0 \\
0 & 1  & -1 & 1 & -1 & 1 & 0 \\
0 & 0  &  1 & -1 & 1 & 0 & 0 \\
0 & 0 & 0 & 1 & 0 & 0 & 0
\end{array}
\right)
$$
is an alternating sign matrix. The fascinating story of alternating sign matrices~\cite{bressoud}
began when combinatorialists where confronted with
a conjecture by Mills, Robbins and Rumsey~\cite{mills2,mills}, which states that the
number of $n \times n$ alternating sign matrices is given by the following simple formula
$$
\prod_{j=0}^{n-1} \frac{(3j+1)!}{(n+j)!}.
$$
For a long time no one could explain this, until finally Zeilberger~\cite{zeilberger} came up with the first
proof. Soon after another, shorter, proof was given by Kuperberg~\cite{kuperberg}. See also \cite{fischer1}, 
where we have recently presented a new proof of this result.

A {\it monotone triangle} is a triangular array $(a_{i,j})_{1 \le j \le i \le n}$ of integers,
$$
\begin{array}{ccccccccccccc}
  &   &            &         &         &         & a_{1,1} &        &         &         &            &     & \\
   &   &           &         &         & a_{2,1} &         &a_{2,2} &         &         &             &    & \\
    &   &           &         & a_{3,1} &         & a_{3,2} &        & a_{3,3} &         &             &    &  \\
     &   &          & a_{4,1} &         & a_{4,2} &         & a_{4,3}&         & a_{4,4} &              &    & \\
      &   & a_{5,1} &         & a_{5,2} &         & a_{5,3} &        & a_{5,4} &         & a_{5,5}       &   &  \\
   & a_{6,1} &      & a_{6,2} &         & a_{6,3} &         & a_{6,4} &        & a_{6,5} &               & a_{6,6} & \\
 a_{7,1} &  & a_{7,2} &    & a_{7,3} &   & a_{7,4} &    & a_{7,5} &     & a_{7,6} & & a_{7,7}
\end{array}
$$
such that $a_{i,j} \le a_{i-1,j} \le a_{i,j+1}$ and $a_{i,j} < a_{i,j+1}$ for all $i,j$. For instance,
$$
\begin{array}{ccccccccccccc}
  &   &            &         &         &         & 4 &        &         &         &            &     & \\
   &   &           &         &         & 2 &         & 6 &         &         &             &    & \\
    &   &           &         & 1 &         & 4 &        & 7 &         &             &    &  \\
     &   &          & 1 &         & 3 &         & 5&         & 7 &              &    & \\
      &   & 1 &         & 2 &         & 4 &        & 6 &         & 7       &   &  \\
   & 1 &      & 2 &         & 3 &         & 5 &        & 6 &               & 7 & \\
 1 &  & 2 &    & 3 &   & 4 &    & 5 &     & 6 & & 7
\end{array}
$$
is a monotone triangle. It corresponds to the alternating sign
matrix above: in the matrix, replace every entry with the sum of
entries in the same column above, the entry itself included. The
result is a $0$-$1$-matrix with one $1$ in the first row, two $1$s
in the second row etc. If one records the columns of the $1$s
rowwise, one obtains the corresponding monotone triangle. It is
not hard to see that this establishes a bijection between monotone
triangles with bottom row $(1,2,\ldots,n)$ and $n \times n$
alternating sign matrices.

Observe that the alternating sign matrix given above is symmetric with respect to the vertical symmetry
axis. This is not the case for all alternating sign matrices. In fact there only exist vertically
symmetric alternating sign matrices of odd size. (This follows from the fact that an alternating 
sign matrix has always a unique $1$ in its top row.) Kuperberg~\cite{kuperberg2} showed that the number of
vertically symmetric $(2n+1) \times (2n+1)$ alternating sign matrices is given by
$$
\frac{n!}{(2n)! 2^n} \prod_{j=1}^n \frac{(6j-2)!}{(2n+2j-1)!}.
$$
(This formula was conjectured by Robbins~\cite{robbins2}.)
The symmetry of vertically symmetric alternating sign matrices translates into a symmetry of the
corresponding monotone triangle: the replacement of every entries $e$ by $2n+2-e$ and the subsequent reflection along the vertical symmetry axis leaves the
monotone triangle invariant. Thus, in case of a vertically symmetric alternating sign matrix,
it suffices to ``store'' (a bit less than) half of the monotone triangle. In our example, this is the following
array.
$$
\begin{array}{cccccc}
         &         &         &         &         & 2 \\
         &         &         &         & 1 &         \\
         &         &         & 1 &         & 3 \\
         &         & 1 &         & 2 &         \\
         & 1 &         & 2 &         & 3 \\
 1 &         & 2 &         & 3 &
\end{array}
$$
Note that the middle column of a monotone triangle corresponding to a
$(2n+1) \times (2n+1)$ vertically symmetric alternating sign matrix consists
solely of $(n+1)$s and, consequently, we do not have to store it.

These considerations led us to the following definition. A {\it halved monotone triangle} is
a triangular array $(a_{i,j})_{1 \le i \le n, 1 \le j \le \lceil \frac{i}{2} \rceil}$ of integers,
$$
\begin{array}{cccccc}
         &         &         &         &         & a_{1,1} \\
         &         &         &         & a_{2,1} &         \\
         &         &         & a_{3,1} &         & a_{3,2} \\
         &         & a_{4,1} &         & a_{4,2} &         \\
         & a_{5,1} &         & a_{5,2} &         & a_{5,3} \\
 a_{6,1} &         & a_{6,2} &         & a_{6,3} &
\end{array}
$$
which is monotone increasing in northeast and southeast direction and strictly increasing along rows, that is
$a_{i+1,j} \le a_{i,j}$,  $a_{i,j} \le a_{i+1,j+1}$ and  $a_{i,j} < a_{i,j+1}$ for all $i,j$. The bijection
sketched above shows that halved monotone triangles
$(a_{i,j})_{1 \le i \le 2n, 1 \le j \le \lceil \frac{i}{2} \rceil}$ with bottom row $(1,2,\ldots,n)$
such that no entry is greater than $n$ correspond to $(2n+1) \times (2n+1)$ vertically symmetric alternating
sign matrices. We are ready to state the main result of the paper.

\begin{theo}
\label{main}
The number of halved monotone triangles with $n$ rows, where no entry exceeds $x$ and with 
bottom row $(k_1,\ldots,k_{\lceil n/2 \rceil})$, $k_1 < k_2 < \ldots < k_{\lceil n/2 \rceil} \le x$, 
is equal to
\begin{multline*}
\left( \prod_{1 \le p < q \le (n+1)/2}
E_{k_p} (E^{-1}_{k_p} + E^{-1}_{k_q} - \id) (E^{-1}_{k_p} + E_{k_q} - \id) \right) \\
\prod_{1 \le i < j \le (n+1)/2} \frac{(k_j - k_i) (2x+1-k_i-k_j)}{(j-i)(j+i-1)}
\end{multline*}
if $n$ is odd and  equal to
\begin{multline*}
\left( \prod_{1 \le p < q \le n/2}
E_{k_p} (E^{-1}_{k_p} + E^{-1}_{k_q} - \id) (E^{-1}_{k_p} + E_{k_q} - \id) \right) \\
\prod_{1 \le i < j \le n/2} \frac{(k_j - k_i) (2x+2-k_i-k_j)}{(j-i)(j+i)}  \prod_{i=1}^{n/2} \frac{x+1-k_i}{i}
\end{multline*}
if $n$ is even, where $E_x p(x) = p(x+1)$ denotes the shift operator. In this formula, 
the product of operators is 
the composition, and, since the shift operators with 
respect to different variables commute, we do not have to specify the order in which they are applied.
\end{theo}

\medskip

We think that the following phenomenon is interesting, since it is in analogy to the
situation for ordinary monotone triangles, see \cite{fischer}. If we consider
``halved monotone triangles'' which do not necessarily have strict monotony 
along rows (the weak monotony follows from the other conditions), the enumeration 
problem is significantly easier. These objects are 
equivalent to shifted plane partitions of trapezoidal shape with
prescribed diagonal and were enumerated by Proctor~\cite[Prop. 4.1]{proctor}. 
The number of these halved triangles with
$n$ rows, where no entry exceeds $x$ and with bottom row $(k_1,\ldots,k_{\lceil n/2 \rceil})$ 
is equal to
$$
\prod_{1 \le i < j \le
(n+1)/2} \frac{(k_j - k_i + j - i) (2x+2+n-i-j-k_i-k_j)}{(j-i)(j+i-1)}
$$
if $n$ is odd and equal to
$$
\prod_{1 \le i < j \le
n/2} \frac{(k_j-k_i+j-i)(2x+2-i-j+n-k_i-k_j)}{(j-i)(j+i)} \prod_{i=1}^{n/2}
\frac{x+1-i+n/2-k_i}{i}
$$
if $n$ is even. Let $\beta(n,x;k_1,\ldots,k_{\lceil n/2 \rceil})$
denote the number of these objects. Then, by Theorem~\ref{main},
the number of halved monotone triangles with $n$ rows, where no entry 
exceeds $x$ and with bottom row $(k_1,\ldots,k_{\lceil n/2 \rceil})$ 
is given by
$$
\left( \prod_{1 \le p < q \le \lceil n/2 \rceil }  (E_{k_p} + E_{k_q} - E_{k_p} E_{k_q}) 
(E_{k_p} + E^{-1}_{k_q} - E_{k_p} E^{-1}_{k_q})  \right)
\beta(n,x;k_1,\ldots,k_{\lceil n/2 \rceil}).
$$
This happens to be in perfect analogy to the situation for 
ordinary monotone triangles: an enumeration formula for the objects with 
strict monotony along rows can be obtained by applying a product of simple 
operators to an enumeration formula for the corresponding objects with weak monotony along rows, the
latter of which is a simple product formula. 

\medskip

The paper is organized as follows.  We prove 
Theorem~\ref{main} in Sections~\ref{rec} -- \ref{derivation}. Our strategy is to first 
show the polynomiality of the formula, then compute its degree, and finally derive enough properties
that characterize the polynomial. To be more 
precise, in Section~\ref{rec}, we introduce the recursion underlying our enumeration formula
for halved monotone triangles and show the polynomiality of the formula. In Section~\ref{operator}, 
we define an operator, which is closely related to the 
recursion and prove a number of lemmas on it. In Section~\ref{properties}, we list and derive the properties 
that characterize the 
enumeration polynomial, and, in Section~\ref{derivation}, we finally 
use these properties
to prove Theorem~\ref{main}. In Section~\ref{generating}, we use our operator formulas 
to derive a generating function for halved monotone triangles as 
well as a generating function for ordinary monotone triangles. 

\section{A recursion and the polynomiality of the enumeration formula}
\label{rec}

For $n \ge 1$ and $k_1 < k_2 < \ldots < k_{\lceil \frac{n}{2} \rceil} \le x$, let
$\gamma(n,x;k_1,\ldots,k_{\lceil \frac{n}{2} \rceil})$ denote the quantity we want to compute, i.e. the 
number of halved monotone triangles
with $n$ rows, where the bottom row is $(k_1,k_2,\ldots,k_{\lceil \frac{n}{2} \rceil})$ and all entries are
no greater than $x$. We define a summation operator for functions $f(l_1,\ldots,l_{m-1})$,  where $m \ge 2$
and $(l_1,\ldots,l_{m-1}) \in \mathbb{Z}^{m-1}$,  as follows. For given 
$(k_1,\ldots,k_m) \in \mathbb{Z}^m$ we have 
\begin{equation}
\label{sumop} \sum^{(k_1,\ldots,k_m)}_{(l_1,\ldots,l_{m-1})} :=
\sum_{(l_1,\ldots,l_{m-1}) \in \mathbb{Z}^{m-1}, \atop  k_1 \le
l_1 \le k_2 \le \ldots \le k_{m-1} \le l_{m-1} \le k_{m}, l_i
\not= l_{i+1}}, 
\end{equation}
i.e. we sum over all strictly increasing sequences $(l_1,\ldots,l_{m-1})$ such that 
$k_i \le l_i \le k_{i+1}$ for all $i$.
This operator is well-defined for all strictly increasing sequences $(k_1,\ldots,k_m) \in \mathbb{Z}^m$.
If we define $\gamma(0,x;-)=1$ we have the following recursions. If $n$ is even then
$$
\gamma(n,x;k_1,\ldots, k_{n/2}) =
\sum_{ (l_1,\ldots, l_{n/2})}^{(k_1,\ldots,k_{n/2},x)}
\gamma(n-1,x;l_1,l_2,\ldots,l_{n/2})
$$
and if $n$ is odd then
$$
\gamma(n,x;k_1,\ldots, k_{(n+1)/2}) =
\sum_{ (l_1,\ldots, l_{(n-1)/2})}^{(k_1,\ldots,k_{(n+1)/2})}
 \gamma(n-1,x;l_1,\ldots,l_{(n-1)/2}).
$$
We want to extend the interpretation of
$\gamma(n,x;k_1,\ldots,k_{\lceil \frac{n}{2} \rceil})$ to
arbitrary $(k_1,\ldots,k_{\lceil \frac{n}{2} \rceil})  \in
\mathbb{Z}^{\lceil \frac{n}{2} \rceil}$. For this purpose, it
suffices to extend the definition of \eqref{sumop} to arbitrary
$(k_1,\ldots,k_m) \in \mathbb{Z}^m$ and then use the recursions to
define the generalization of $\gamma$. We use induction with respect to $m$. For $m=2$,
let
$$
\sum_{(l_1)}^{(k_1,k_2)} a(l_1) := \sum_{l_1=k_1}^{k_2} a(l_1),
$$
where here and in the following $\sum\limits_{i=a}^b f(i) = - \sum\limits_{i=b+1}^{a-1} f(i)$ if $a > b$. (Note that
this implies $\sum\limits_{i=a}^{a-1} f(i) = 0$. Moreover, $\sum\limits_{x=0}^y p(x)$ will be a polynomial in $y$ if $p(x)$ is
a polynomial in $x$.) If $m > 2$ we define
\begin{multline*}
\sum_{(l_1,\ldots,l_{m-1})}^{(k_1,\ldots,k_m)} a(l_1,\ldots,l_{m-1}) := \\
\sum_{(l_1,\ldots,l_{m-2})}^{(k_1,\ldots,k_{m-1})} \sum_{l_{m-1}=k_{m-1}+1}^{k_m} a(l_1,\ldots, l_{m-2}, l_{m-1}) +
\sum_{(l_1,\ldots,l_{m-2})}^{(k_1,\ldots,k_{m-2}, k_{m-1}-1)} a(l_1,\ldots,l_{m-2},k_{m-1}).
\end{multline*}
Now it is also obvious that $\gamma(n,x;k_1,\ldots, k_{\lceil n/2 \rceil})$ is a polynomial in
$(k_1,\ldots, k_{\lceil n/2 \rceil},x)$ for fixed $n$.

This recursion can be used to compute $\gamma(n,x;k_1,\ldots,k_{\lceil n/2 \rceil})$ for 
small values of $n$. For $n=1,2,3,4,5$ you find the results below.

\begin{multline*} 1, x-k_{1}+1, \frac{1}{2} (2 x+2-k_{1}-k_{2}) (k_{2}-k_{1}+1), 
\frac{1}{6} (x-k_2+1) (-k_{1}^3+3 x k_{1}^2 +6 k_{1}^2 -2 x^2 k_{1}  \\ +k_{2}^2 k_{1}  -10 x
   k_{1}  -2 x k_{2} k_{1}-2 k_{2} k_{1}-11 k_{1}+2 x^2-x k_{2}^2-2 k_{2}^2+7 x+2 x^2 k_{2} 
 +6 x  k_{2}+4 k_{2} ), \\
\frac{1}{48} (k_{2}^2 k_{1}^4-k_{3}^2 k_{1}^4+2 x k_{1}^4-2 x k_{2} k_{1}^4-3 k_{2} k_{1}^4+2
   x k_{3} k_{1}^4  +k_{3} k_{1}^4+2 k_{1}^4-8 x^2 k_{1}^3-4 x k_{2}^2 k_{1}^3 \\  -10 k_{2}^2 k_{1}^3+4 x
   k_{3}^2 k_{1}^3 +10 k_{3}^2 k_{1}^3-28 x k_{1}^3+8 x^2 k_{2} k_{1}^3  +32 x k_{2} k_{1}^3+30 k_{2}
   k_{1}^3-8 x^2 k_{3} k_{1}^3-24 x k_{3} k_{1}^3\\ -10 k_{3} k_{1}^3-20 k_{1}^3  -k_{2}^4
   k_{1}^2+k_{3}^4 k_{1}^2  +8 x^3 k_{1}^2+4 x k_{2}^3 k_{1}^2  +6 k_{2}^3 k_{1}^2-4 x k_{3}^3
   k_{1}^2-2 k_{3}^3 k_{1}^2+60 x^2 k_{1}^2+12 x k_{2}^2 k_{1}^2 \\ +24 k_{2}^2 k_{1}^2  -24 x k_{3}^2
   k_{1}^2  -36 k_{3}^2 k_{1}^2+122 x k_{1}^2-8 x^3 k_{2} k_{1}^2-60 x^2 k_{2} k_{1}^2-138 x k_{2}
   k_{1}^2  -99 k_{2} k_{1}^2+8 x^3 k_{3} k_{1}^2 \\ +60 x^2 k_{3} k_{1}^2+102 x k_{3} k_{1}^2  +37 k_{3}
   k_{1}^2+70 k_{1}^2+2 x k_{2}^4 k_{1}+5 k_{2}^4 k_{1}-2 x k_{3}^4 k_{1}-5 k_{3}^4 k_{1}-24 x^3
   k_{1}-8 x^2 k_{2}^3 k_{1} \\ -32 x k_{2}^3 k_{1}  -30 k_{2}^3 k_{1}+8 x^2 k_{3}^3 k_{1}  +24 x k_{3}^3
   k_{1}+10 k_{3}^3 k_{1}-124 x^2 k_{1}+8 x^3 k_{2}^2 k_{1}+36 x^2 k_{2}^2 k_{1}+42 x k_{2}^2 
   k_{1} \\ +29 k_{2}^2 k_{1}-8 x^3 k_{3}^2 k_{1} -12 x^2 k_{3}^2 k_{1}  +42 x k_{3}^2 k_{1}+31 k_{3}^2
   k_{1}-152 x k_{1}+16 x^3 k_{2} k_{1}+96 x^2 k_{2} k_{1}+140 x k_{2} k_{1}\\ +48 k_{2} k_{1}-32 x^3
   k_{3} k_{1}-144 x^2 k_{3} k_{1}  -136 x k_{3} k_{1}-36 k_{3} k_{1}-52 k_{1}-4 x k_{2}^4-6
   k_{2}^4-k_{2}^2 k_{3}^4+2 x k_{3}^4\\ +2 x k_{2} k_{3}^4  +3 k_{2} k_{3}^4 +4 k_{3}^4+16 x^2
   k_{2}^3+48 x k_{2}^3+28 k_{2}^3  -8 x^2 k_{3}^3 + 4 x k_{2}^2 k_{3}^3  +2 k_{2}^2 k_{3}^3-20 x
   k_{3}^3-8 x^2 k_{2} k_{3}^3 \\ -16 x k_{2} k_{3}^3  -6 k_{2} k_{3}^3-8 k_{3}^3-16 x^3 k_{2}^2-72 x^2
   k_{2}^2-100 x k_{2}^2  -66 k_{2}^2+k_{2}^4 k_{3}^2+8 x^3 k_{3}^2 -4 x k_{2}^3 k_{3}^2-6 k_{2}^3
   k_{3}^2 \\ +12 x^2 k_{3}^2 +12 x k_{2}^2 k_{3}^2 +12 k_{2}^2 k_{3}^2-22 x k_{3}^2+8 x^3 k_{2}
   k_{3}^2+12 x^2 k_{2} k_{3}^2  -6 x k_{2} k_{3}^2  +15 k_{2} k_{3}^2-4 k_{3}^2+32 x^2 k_{2} \\ +96 x
   k_{2}+44 k_{2}-2 x k_{2}^4 k_{3} -k_{2}^4 k_{3}+24 x^3 k_{3}+8 x^2 k_{2}^3 k_{3} +16 x k_{2}^3
   k_{3}+6 k_{2}^3 k_{3}  +92 x^2 k_{3}-8 x^3 k_{2}^2 k_{3} \\ -36 x^2 k_{2}^2 k_{3} -42 x k_{2}^2
   k_{3}-13 k_{2}^2 k_{3}+56 x k_{3}+16 x^3 k_{2} k_{3}+48 x^2 k_{2} k_{3}-4 x k_{2} k_{3}  -12 k_{2}
   k_{3}+8 k_{3}) 
\end{multline*}
This data suggests that the degree of 
$\gamma(n,x;k_1,\ldots,k_{\lceil n/2})$ in $k_i$ is always $n-1$. That this is indeed the case will be 
shown in the following section. However, this comes by surprise because of the following: suppose 
that $a(l_1,\ldots,l_{m-1})$ is a polynomial of degree $R$ in every $l_i$. Then the degree
of 
$$
\sum_{(l_1,\ldots,l_{m-1})}^{(k_1,\ldots,k_m)}  a(l_1,\ldots,l_{m-1})
$$
in $k_i$ could be as high as $2 R + 2$ (e.g. $a(l_1,\ldots,l_{m-1})=\prod\limits_{i=1}^{m-1} l_i^R$). 
This estimation provides (by induction with respect to $n$) a bound of 
$2^{n+1}-2$ for the degree of $\gamma(n,x;k_1,\ldots,k_{\lceil n/2 \rceil})$ in $k_i$.

\section{An operator related to the recursion}
\label{operator}

Most of the definitions and lemmas in this section are taken from \cite{fischer}.  The {\it shift
operator} $E_x$ is defined as $E_x p(x) = p(x+1)$ and the {\it difference operator} $\Delta_x$ is
defined as $E_x - \id$. The {\it swapping operator} $S_{x,y}$ is defined as $S_{x,y} f(x,y) = f(y,x)$.

Note that the shift operator $E_x$ is invertible as an operator over $\mathbb{C}[x]$, whereas the difference operator
$\Delta_x$ is not, since it decreases the degree of a polynomial. In the following, we will
consider rational functions in shift operators and thus we need a lemma in order to show that the inverses of our denominators exist. (The lemma is a generalization of \cite[Lemma~1]{fischer}.) For the statement of
the lemma we need the following observation. Let
$$
p(X_1,\ldots,X_n) = \sum_{(i_1,\ldots,i_n) \in (\mathbb{Z}^{\ge 0})^n} a_{i_1,\ldots,i_n} X_1^{i_1}  \cdots X_n^{i_n}
$$
be a formal power series in $(X_1,\ldots,X_n)$ over $\mathbb{C}$ and $G(k_1,\ldots,k_n)$ be a polynomial
in $(k_1,\ldots,k_n)$ over $\mathbb{C}$. We define
$$
p(\Delta_{k_1},\ldots,\Delta_{k_n}) G(k_1,\ldots,k_n):  = \\
\sum_{(i_1,\ldots,i_n) \in (\mathbb{Z}^{\ge 0})^{n}} a_{i_1,\ldots,i_n} \Delta_{k_1}^{i_1} \cdots \Delta_{k_n}^{i_n}
G(k_1,\ldots,k_n).
$$
This is a finite sum and thus well-defined since $\Delta_{k_i}^{d+1} G(k_1,\ldots,k_n) = 0$ if
$\deg_{k_i} G(k_1,\ldots,k_n) = d$.

\begin{lem}
\label{invers}
Let $p(X_1,\ldots,X_n)$ be a formal power series in $(X_1,\ldots,X_n)$ over $\mathbb{C}$ with
non-zero constant term. Then
$p(\Delta_{k_1},\ldots,\Delta_{k_n})$ is invertible as an operator over $\mathbb{C}[k_1,\ldots,k_n]$, i.e. there
exists a formal power series $q(X_1,\ldots,X_n)$ with 
\begin{multline*}
p(\Delta_{k_1},\ldots,\Delta_{k_n}) q(\Delta_{k_1},\ldots,\Delta_{k_n}) F(k_1,\ldots,k_n) \\ =
q(\Delta_{k_1},\ldots,\Delta_{k_n}) p(\Delta_{k_1},\ldots,\Delta_{k_n}) F(k_1,\ldots,k_n) =
F(k_1,\ldots,k_n)
\end{multline*}
for all polynomials $F(k_1,\ldots,k_n)$. 
Moreover
\begin{multline*}
\deg_{k_{i_1},\ldots,k_{i_m}} G(k_1,\ldots,k_n) = \deg_{k_{i_1},\ldots,k_{i_m}} p(\Delta_{k_1},\ldots,\Delta_{k_n}) G(k_1,\ldots,k_n) =
\\
\deg_{k_{i_1},\ldots,k_{i_m}} q(\Delta_{k_1},\ldots,\Delta_{k_n}) G(k_1,\ldots,k_n)
\end{multline*}
for all $(i_1, i_2, \ldots,i_m)$ with $1 \le i_1 < i_2 < \ldots < i_m \le n$, where $\deg_{k_{i_1},\ldots,k_{i_m}} G(k_1,\ldots,k_n)$ denotes
the degree of $G(k_1,\ldots,k_n)$ as a polynomial in $(k_{i_1},\ldots,k_{i_m})$.
\end{lem}

{\it Proof.} The assertion follows from the fact that $p(X_1,\ldots,X_n)$ is invertible in the (commutative) 
algebra of formal power series over $\mathbb{C}$ if (and only if) $p(X_1,\ldots,X_n)$ has a non-zero constant term.
This is because $p(X_1,\ldots,X_n) q(X_1,\ldots, X_n)=1$ is equivalent to $a_{0,\ldots,0} b_{0,\ldots,0} = 1$ and
$$
\sum_{(i_1,\ldots,i_n), (j_1,\ldots,j_n) \in (\mathbb{Z}^{\ge 0})^n \atop
(i_1,\ldots,i_n) + (j_1,\ldots, j_n) = (r_1,\ldots,r_n)} a_{i_1,\ldots,i_n} b_{j_1,\ldots,j_n} = 0
$$
for $(r_1,\ldots,r_n) \not= (0,\ldots,0)$, where
$$
p(X_1,\ldots,X_n) = \sum_{(j_1,\ldots,j_n) \in (\mathbb{Z}^{\ge 0})^n} a_{j_1,\ldots,j_n} X_1^{j_1}  \cdots X_n^{j_n}.
$$
and
$$
q(X_1,\ldots,X_n) = \sum_{(j_1,\ldots,j_n) \in (\mathbb{Z}^{\ge 0})^n} b_{j_1,\ldots,j_n} X_1^{j_1}  \cdots X_n^{j_n}.
$$
By assumption $a_{0,\ldots,0} \not= 0$ and, consequently, the equations allow us to determine the coefficients
$b_{r_1,\ldots,r_n}$ by induction with respect to $r_1 + \ldots + r_n$. The assertion about the degree follows from
the fact that
\begin{multline*}
\deg_{k_{i_1},\ldots,k_{i_m}} p(\Delta_{k_1},\ldots,\Delta_{k_n}) G(k_1,\ldots,k_n) \le \deg_{k_{i_1},\ldots,k_{i_m}} G(k_1,\ldots,k_n) \\
= \deg_{k_{i_1},\ldots,k_{i_m}} q(\Delta_{k_1},\ldots,\Delta_{k_n}) p(\Delta_{k_1},\ldots,\Delta_{k_n})
G(k_1,\ldots,k_n) \\ \le \deg_{k_{i_1},\ldots,k_{i_m}} p(\Delta_{k_1},\ldots,\Delta_{k_n}) G(k_1,\ldots,k_n). \qed
\end{multline*}

\medskip

We define $V_{x,y}= (\id + E_y \Delta_x) = E_x + \Delta_x \Delta_y$.
In  the following lemma we will see why this operator is
of significance for the recursion underlying $\gamma(n,x;k_1,\ldots,k_{\lceil n/2 \rceil})$. (The lemma is
equivalent to \cite[Lemma~2]{fischer}.) It will be used for showing that the degree of
$\gamma(n,x;k_1,\ldots, k_{\lceil n/2 \rceil})$ is no greater than $n-1$ in every $k_i$.

\begin{lem}
\label{degree}
Let $a(x,y)$ be a polynomial in $x$ and $y$ which is of degree no greater than $R$ in each of $x$ and $y$.
Moreover, assume that $(\id + S_{x,y}) V_{x,y} a(x,y)$ is of degree no greater than $R$ as a polynomial
in $x$ and $y$, i.e. a linear combination of monomials $x^m y^n$ with $m+n \le R$. Then
\begin{equation}
\label{rtf}
\sum_{(x,y)}^{(k_1,k_2,k_3)} a(x,y) = \sum_{x=k_1}^{k_2} \sum_{y=k_2}^{k_3} a(x,y) - a(k_2,k_2)
\end{equation}
is of degree no greater than $R+2$ in $k_2$. Moreover, if $(\id + S_{x,y}) V_{x,y} a(x,y)=0$ then
the degree of \eqref{rtf} in $k_2$ is no greater than $R+1$.
\end{lem}

{\it Proof.} First note that, by Lemma~\ref{invers}, $V_{x,y} + V_{y,x} = 2 \id + \Delta_x + \Delta_y + 2 \Delta_x \Delta_y$
is invertible. Thus
\begin{equation}
\label{rtf1}
(\id + S_{x,y}) \frac{V_{x,y}}{V_{x,y} + V_{y,x}} + \frac{V_{y,x}}{V_{x,y} + V_{y,x}} (\id - S_{x,y}) = \id,
\end{equation}
since $V_{x,y} + V_{y,x}$ and $S_{x,y}$ commute. Moreover,
$$
(\id + S_{x,y}) \frac{V_{x,y}}{V_{x,y} + V_{y,x}} a(x,y) = \frac{1}{V_{x,y}+V_{y,x}} (\id + S_{x,y}) V_{x,y} a(x,y).
$$
By Lemma~\ref{invers}, the degree of this expression in $x$ and $y$ is equal to the
degree of $(\id + S_{x,y}) V_{x,y} a(x,y)$ in $x$ and $y$ and the expression vanishes if and only if
$(\id + S_{x,y}) V_{x,y} a(x,y)$ vanishes. Thus, by \eqref{rtf1}, it suffices to show that the degree of
$$
\sum_{(x,y)}^{(k_1,k_2,k_3)}  \frac{V_{y,x}}{V_{x,y} + V_{y,x}} (\id - S_{x,y}) a(x,y)
$$
in $k_2$ is no greater than $R+1$. Once more the degree estimation from Lemma~\ref{invers} implies that this can
be reduced to showing the following.
If we define $b(x,y) = V_{y,x} (\id - S_{x,y}) \binom{x}{p} \binom{y}{q}$ then the degree of
$
\sum\limits_{(x,y)}^{(k_1,k_2,k_3)} b(x,y)
$
in $k_2$ is no greater than $\max(p,q)+1$. In order to do so, observe that
\begin{multline*}
b(x,y) = V_{y,x} (\id - S_{x,y}) \binom{x}{p} \binom{y}{q} \\
= \binom{x}{p} \binom{y}{q} - \binom{x}{q} \binom{y}{p} + \binom{x+1}{p} \binom{y}{q-1} - \binom{x+1}{q} \binom{y}{p-1}.
\end{multline*}
Therefore, and by the summation formula
$$
\sum_{x=a}^{b} \binom{x}{n} = \sum_{x=a}^b \binom{x+1}{n+1} - \binom{x}{n+1} = \binom{b+1}{n+1} - \binom{a}{n+1},
$$
we have
\begin{multline*}
\sum_{x=k_1}^{k_2} \sum_{y=k_2}^{k_3} b(x,y)  - b(k_2,k_2)
= \left( \binom{k_2+1}{p+1} - \binom{k_1}{p+1} \right) \left( \binom{k_3+1}{q+1} - \binom{k_2}{q+1} \right)  \\
- \left( \binom{k_2+1}{q+1} - \binom{k_1}{q+1} \right) \left( \binom{k_3+1}{p+1} - \binom{k_2}{p+1} \right)  \\
+ \left( \binom{k_2+2}{p+1} - \binom{k_1+1}{p+1} \right) \left( \binom{k_3+1}{q} - \binom{k_2}{q} \right) \\
- \left( \binom{k_2+2}{q+1} - \binom{k_1+1}{q+1} \right) \left( \binom{k_3+1}{p} - \binom{k_2}{p} \right)   \\
- \binom{k_2+1}{p} \binom{k_2}{q-1} + \binom{k_2+1}{q} \binom{k_2}{p-1}.
\end{multline*}
If we repeatedly apply the identity
$$
\binom{n}{k} = \binom{n-1}{k} + \binom{n-1}{k-1}
$$
to this expression, we finally see that this is a polynomial in $k_2$ of degree no greater than $\max(p,q)+1$. \qed

\medskip

In order to use this lemma to compute the degree of $\gamma(n,x;k_1,\ldots,k_{\lceil n/2 \rceil})$
in every $k_i$, we need to show that
$
(\id + S_{k_i,k_{i+1}}) V_{k_i,k_{i+1}} \gamma(n,x;k_1,\ldots,k_{\lceil n/2 \rceil}) = 0
$
for all $i$. This will be a consequence of the following lemma, which implies that
$$
(\id + S_{k_i,k_{i+1}}) V_{k_i,k_{i+1}}
\left( \sum_{(l_1,\ldots,l_{m-1})}^{(k_1,\ldots,k_m)} a(l_1,\ldots,l_{m-1}) \right)
$$
is expressible as a certain sum of $$(\id + S_{l_{i-1},l_i}) V_{l_{i-1},l_i} a(l_1,\ldots,l_{m-1})$$ and
$$(\id + S_{l_{i},l_{i+1}}) V_{l_{i},l_{i+1}} a(l_1,\ldots,l_{m-1}).$$ It is yet another result, 
which manifests the connection of $V_{x,y}$ and the recursion. The lemma is
Lemma~3 of \cite{fischer} and we omit its proof here. In order to simplify the
statement we use the following notation $T_{x,y} = (\id + S_{x,y}) V_{x,y}$.

\begin{lem}
\label{kompliziert}
Let $f(l_1,l_2,l_3)$ be a function on $\mathbb{Z}^3$ with values in
$\mathbb{C}$ and define
$$
g(k_1,k_2,k_3,k_4)= \sum_{(l_1,l_2,l_3)}^{(k_1,k_2,k_3,k_4)} f(l_1,l_2,l_3).
$$
Then
\begin{multline*}
T_{k_2,k_3} g(k_1,k_2,k_3,k_4)  \\
= - \frac{1}{2} \left(
\sum_{l_1=k_2+1}^{k_3} \sum_{l_2=k_2+1}^{k_3} \sum_{l_3=k_2}^{k_4}
T_{l_1,l_2} f(l_1,l_2,l_3) +
\sum_{l_1=k_1}^{k_2+1} \sum_{l_2=k_2}^{k_3-1} \sum_{l_3=k_2}^{k_3-1}
T_{l_2,l_3} f(l_1,l_2,l_3) \right) \\
+ \frac{1}{2}
\left( \sum_{l_1=k_2}^{k_3-1} \sum_{l_2=k_2}^{k_3-1} \Delta_{l_2} (\id + E_{l_1})
T_{l_1,l_2} f(l_1,l_2,k_2) -
\sum_{l_2=k_2}^{k_3-1} \sum_{l_3=k_2}^{k_3-1} \Delta_{l_2} (\id + E_{l_3} ) T_{l_2,l_3} f(k_2+1,l_2,l_3) \right) \\
+ \frac{1}{2} \left(
\left. T_{l_1,l_2} f(l_1,l_2,k_2+1) \right|_{(l_1,l_2)=(k_2,k_2)}
- \left. T_{l_1,l_2} f(l_1,l_2,k_3+1) \right|_{(l_1,l_2)=(k_2,k_2)} \right. \\
+ \left. \left. T_{l_2,l_3} f(k_2,l_2,l_3) \right|_{(l_2,l_3)=(k_2,k_2)}
- \left. T_{l_2,l_3} f(k_3,l_2,l_3) \right|_{(l_2,l_3)=(k_2,k_2)} \right) \\
-  \left. T_{l_1,l_2} f(l_1,l_2,k_2+1) \right|_{(l_1,l_2)=(k_2,k_3)}
-  \left. T_{l_2,l_3} f(k_2,l_2,l_3) \right|_{(l_2,l_3)=(k_2,k_3)}.
\end{multline*}
Moreover, for a function $h(l_1,l_2)$ on $\mathbb{Z}^2$,
$$
T_{k_1,k_2} \sum_{(l_1,l_2)}^{(k_1,k_2,k_3)} h(l_1,l_2) =
- \frac{1}{2} \sum_{l_1=k_1}^{k_2-1} \sum_{l_2=k_1}^{k_2-1} T_{l_1,l_2} h(l_1,l_2).
$$
\end{lem}

This proves the assertion preceding the lemma for $m=3,4$. For $m=2$ observe that
\begin{multline*}
(\id + S_{k_1,k_2}) V_{k_1,k_2} \sum_{(l_1)}^{(k_1,k_2)} a(l_1)
= (\id + S_{k_1,k_2}) \left( \sum_{l_1=k_1}^{k_2} a(l_1) - a(k_1)\right) \\
= \sum_{l_1=k_1+1}^{k_2} a(l_1) + \sum_{l_2=k_2+1}^{k_1} a(l_1) = 0.
\end{multline*}
In order to use Lemma~\ref{kompliziert} to prove the assertion for $m \ge 5$, we need
a merging rule for \eqref{sumop}. Let $f(x,z)$ be a function on $\mathbb{Z}^2$. Then
the operator $I_{x,z}^y$ is defined as follows.
$$
I_{x,z}^y f(x,z) = f(y-1,y) + f(y,y+1) - f(y-1,y+1) = \left. V_{x,z} f(x,z) \right|_{(x=y-1,z=y)}
$$
Using this operator, we have
\begin{equation}
\label{merge}
\sum_{(l_1,\ldots,l_{m-1})}^{(k_1,\ldots,k_m)} a(l_1,\ldots,l_{m-1}) =
I_{w,x}^{k_{i-1}} I_{y,z}^{k_{i+2}}
\sum_{(l_1,\ldots,l_{i-2})}^{(k_1,\ldots,k_{i-2},w)}
\sum_{(l_{i-1},l_i,l_{i+1})}^{(x,k_i,k_{i+1},y)}
\sum_{(l_{i+2},\ldots,l_n)}^{(z,k_{i+3},\ldots,l_{m-1})} a(l_1,\ldots,l_{m-1})
\end{equation}
and this enables one to prove the assertion for $m \ge 5$. (For details see \cite[Section~4]{fischer}.)

After noting that $(\id + S_{x,y}) g(x,y) = 0$ if and only if $g(x,y)$ is antisymmetric in $x$ and $y$, we finally obtain the
following.

\begin{cor}
\label{zero}
Suppose $a(l_1,\ldots,l_{m-1})$ is a function on $\mathbb{Z}^{m-1}$ such that
$V_{l_i,l_{i+1}} a(l_1,\ldots,l_{m-1})$ is antisymmetric in $l_i$ and $l_{i+1}$ for all $i$. Then
$$
V_{k_i,k_{i+1}}  \sum_{(l_1,\ldots,l_{m-1})}^{(k_1,\ldots,k_m)} a(l_1,\ldots,l_{m-1})
$$
is antisymmetric in $k_i$ and $k_{i+1}$ for all $i$.
\end{cor}

\section{Characterizing properties of $\gamma(n;k_1,\ldots,k_{\lceil n/2 \rceil})$}
\label{properties}

We apply the results from the previous section to
$\gamma(n,x;k_1,\ldots,k_{\lceil n/2 \rceil})$: Corollary~\ref{zero} implies 
by induction with respect to $n$ that $V_{k_i,k_{i+1}}
\gamma(n,x;k_1,\ldots,k_{\lceil n/2 \rceil})$ is antisymmetric in
$k_i$ and $k_{i+1}$ for all $i$. Lemma~\ref{degree} and the 
merging rule \eqref{merge} then implies
by induction with respect to $n$ that the degree of
$\gamma(n,x;k_1,\ldots,k_{\lceil n/2 \rceil})$ in $k_i$ is no
greater than $n-1$. We summarize this in the following lemma.

\begin{lem} 
\label{prop1} For $n \ge 1$, $\gamma(n,x;k_1,\ldots,k_{\lceil n/2 \rceil})$ is a polynomial of degree no
greater than $n-1$ in every $k_i$. Furthermore,
$V_{k_i,k_{i+1}} \gamma(n,x;k_1,\ldots,k_{\lceil n/2 \rceil})$ is
antisymmetric in $k_i$ and $k_{i+1}$ for all $i$.
\end{lem}

It will be shown that the properties from the previous lemma together with the property in the following
lemma characterize $\gamma$ up to a multiplicative rational constant.

\begin{lem}
\label{odd-even}
If $n$ is even then
$$
\gamma(n,x;k_1,\ldots,k_{n/2-1},k_{n/2}) =
- \gamma(n,x;k_1,\ldots,k_{n/2-1}, 2x+2-k_{n/2})
$$
and if $n$ is odd then
$$
\gamma(n,x;k_1,\ldots,k_{(n-1)/2}, k_{(n+1)/2}) = \gamma(n,x;k_1,\ldots, k_{(n-1)/2}, 2x+1-k_{(n+1)/2}).
$$
\end{lem}

In order to prove this lemma, we need another lemma.

\begin{lem}
\label{hilf}
(1) Let $f(l_1)$ be such that $f(l_1)=-f(2x+2-l_1)$ for all $l_1$ and 
define $g(k_1,k_2)=\sum\limits_{l_1=k_1}^{k_2} f(l_1)$.
Then $g(k_1,k_2)=g(k_1,2x+1-k_2)$ for all $k_1, k_2$.

(2) Let $f(l_1,l_2)$ be such that $f(l_1,l_2)=f(l_1,2x+1-l_2)$ for all $l_1, l_2$ and
$$(\id + S_{l_1,l_2}) V_{l_1,l_2} f(l_1,l_2)=0.$$
Define
$$
g(k_1,k_2)=\sum_{(l_1,l_2)}^{(k_1,k_2,x)} f(l_1,l_2) = \sum_{l_1=k_1}^{k_2} \sum_{l_2=k_2}^{x} f(l_1,l_2) - f(k_2,k_2).
$$
Then $g(k_1,k_2)=-g(k_1,2x+2-k_2)$ for all $k_1, k_2$.
\end{lem}

{\it Proof of Lemma~\ref{hilf}.} (1) By definition,
$$
g(k_1,2x+1-k_2) = \sum_{l_1=k_1}^{2x+1-k_2} f(l_1) = \sum_{l_1=k_1}^{k_2} f(l_1) + \sum_{l_1=k_2+1}^{2x+1-k_2} f(l_1).
$$
The assertion follows since $\sum\limits_{l_1=k_2+1}^{2x+1-k_2} f(l_1)=0$. This is because
$$
\sum_{l_1=k_2+1}^{2x+1-k_2} f(l_1) = \sum_{l_1=k_2+1}^{2x+1-k_2} - f(2x+2-l_1) = -
\sum_{l_1=k_2+1}^{2x+1-k_2} f(l_1).
$$

(2) Observe that
\begin{multline}
\label{expr}
g(k_1,2x+2-k_2) = \sum_{l_1=k_1}^{2x+2-k_2} \sum_{l_2=2x+2-k_2}^{x} f(l_1,l_2) - f(2x+2-k_2,2x+2-k_2) \\
= \sum_{l_1=k_1}^{k_2} \sum_{l_2=2x+2-k_2}^{x} f(l_1,l_2) + 
\sum_{l_1=k_2+1}^{2x+2-k_2} \sum_{l_2=2x+2-k_2}^{x} f(l_1,l_2) - f(2x+2-k_2,2x+2-k_2) \\
= \sum_{l_1=k_1}^{k_2} \sum_{l_2=2x+2-k_2}^{x} f(l_1,2x+1-l_2) 
-\frac{1}{2} \sum_{l_1=k_2+1}^{2x+2-k_2} \sum_{l_2=x+1}^{2x+1-k_2} f(l_1,l_2)
-\frac{1}{2} \sum_{l_1=k_2+1}^{2x+2-k_2} \sum_{l_2=x+1}^{2x+1-k_2} f(l_1,l_2)
\\ - f(2x+2-k_2,2x+2-k_2) \\
= \sum_{l_1=k_1}^{k_2} \sum_{l_2=x+1}^{k_2-1} f(l_1,l_2) 
-\frac{1}{2} \sum_{l_1=k_2+1}^{2x+2-k_2} \sum_{l_2=x+1}^{2x+1-k_2} f(l_1,l_2)
-\frac{1}{2} \sum_{l_1=k_2+1}^{2x+2-k_2} \sum_{l_2=x+1}^{2x+1-k_2} f(l_1,2x+1-l_2)
\\ - f(2x+2-k_2,2x+2-k_2) \\
= -\sum_{l_1=k_1}^{k_2} \sum_{l_2=k_2}^{x} f(l_1,l_2) 
-\frac{1}{2} \sum_{l_1=k_2+1}^{2x+2-k_2} \sum_{l_2=x+1}^{2x+1-k_2} f(l_1,l_2)
-\frac{1}{2} \sum_{l_1=k_2+1}^{2x+2-k_2} \sum_{l_2=k_2}^{x} f(l_1,l_2)
\\ - f(2x+2-k_2,2x+2-k_2) \\
= -\sum_{l_1=k_1}^{k_2} \sum_{l_2=k_2}^{x} f(l_1,l_2) 
-\frac{1}{2} \sum_{l_1=k_2+1}^{2x+2-k_2} \sum_{l_2=k_2}^{2x+1-k_2} f(l_1,l_2)
- f(2x+2-k_2,2x+2-k_2).
\end{multline}
Moreover, we have
$
E_{x}^{-1} (\id + S_{x,y}) V_{x,y} = (\id + E_y E_x^{-1} S_{x,y}) (\id + E_x^{-1} \Delta_x \Delta_y),
$
and, therefore,
$
(\id + E_{l_2} E_{l_1}^{-1} S_{l_1,l_2}) (\id + E_{l_1}^{-1} \Delta_{l_1} \Delta_{l_2}) f(l_1,l_2) = 0
$.
This implies that
$$
\sum_{l_1=k_2+1}^{2x+2-k_2} \sum_{l_2=k_2}^{2x+1-k_2} (\id + E_{l_1}^{-1} \Delta_{l_1} \Delta_{l_2}) f(l_1,l_2) = 0,
$$
since
\begin{multline*}
\sum_{l_1=k_2+1}^{2x+2-k_2} \sum_{l_2=k_2}^{2x+1-k_2} (\id + E_{l_1}^{-1} \Delta_{l_1} \Delta_{l_2}) f(l_1,l_2) \\
= \sum_{l_1=k_2+1}^{2x+2-k_2} \sum_{l_2=k_2}^{2x+1-k_2} - E_{l_2} E_{l_1}^{-1} S_{l_1,l_2}
(\id + E_{l_1}^{-1} \Delta_{l_1} \Delta_{l_2}) f(l_1,l_2) \\
= - \sum_{l_1=k_2}^{2x+1-k_2} \sum_{l_2=k_2+1}^{2x+2-k_2} S_{l_1,l_2} (\id + E_{l_1}^{-1} \Delta_{l_1} \Delta_{l_2}) f(l_1,l_2) \\
= - \sum_{l_2=k_2}^{2x+1-k_2} \sum_{l_1=k_2+1}^{2x+2-k_2} (\id + E_{l_1}^{-1} \Delta_{l_1} \Delta_{l_2}) f(l_1,l_2).
\end{multline*}
Thus, \eqref{expr} is equal to
\begin{multline}
\label{expr2}
- \sum_{l_1=k_1}^{k_2} \sum_{l_2=k_2}^{x} f(l_1,l_2)+ \frac{1}{2} \sum_{l_1=k_2+1}^{2x+2-k_2} \sum_{l_2=k_2}^{2x+1-k_2} E_{l_1}^{-1} \Delta_{l_1} \Delta_{l_2} f(l_1,l_2)
- f(2x+2-k_2,2x+2-k_2) \\
= - \sum_{l_1=k_1}^{k_2} \sum_{l_2=k_2}^{x} f(l_1,l_2) - \frac{1}{2} f(2x+2-k_2,2x+2-k_2) +
\frac{1}{2} f(k_2,k_2) - \frac{1}{2} f(2x+2-k_2,k_2)  \\ - \frac{1}{2} f(k_2,2x+2-k_2).
\end{multline}
Next observe that
\begin{multline*}
0=\left. \left( (\id + S_{l_1,l_2}) V_{l_1,l_2}  f(l_1,l_2) \right) \right|_{(l_1,l_2)=(k_2-1,2x+1-k_2)} \\ =
f(k_2-1,2x+1-k_2) + f(k_2,2x+2-k_2) - f(k_2-1,2x+2-k_2) \\ + f(2x+1-k_2,k_2-1) + f(2x+2-k_2,k_2) - f(2x+1-k_2,k_2).
\end{multline*}
We replace $- \frac{1}{2} f(2x+2-k_2,k_2) - \frac{1}{2}  f(k_2,2x+2-k_2)$ in \eqref{expr2} by
$$
\frac{1}{2} f(k_2-1,2x+1-k_2) + \frac{1}{2} f(2x+1-k_2,k_2-1) -
\frac{1}{2} f(k_2-1,2x+2-k_2) - \frac{1}{2} f(2x+1-k_2,k_2)
$$
and, consequently, \eqref{expr2} is equal to
\begin{multline}
\label{expr3}
- \sum_{l_1=k_1}^{k_2} \sum_{l_2=k_2}^{x} f(l_1,l_2) + f(k_2,k_2)
- \frac{1}{2} ( f(2x+2-k_2,2x+2-k_2) + f(k_2,k_2) \\
- f(k_2-1,2x+1-k_2) - f(2x+1-k_2,k_2-1) + f(k_2-1,2x+2-k_2) + f(2x+1-k_2,k_2) ) \\
= - \sum_{l_1=k_1}^{k_2} \sum_{l_2=k_2}^{x} f(l_1,l_2) + f(k_2,k_2)
- \frac{1}{2} ( f(2x+2-k_2,2x+2-k_2) + f(k_2,k_2) \\
- f(k_2-1,k_2) - f(2x+1-k_2,2x+2-k_2) + f(k_2-1,k_2-1) + f(2x+1-k_2,2x+1-k_2) ).
\end{multline}
Finally, 
$$
0 = \left. \left( (\id + S_{l_1,l_2}) V_{l_1,l_2}  f(l_1,l_2) \right) \right|_{(l_1,l_2)=(l,l)} =
2 ( f(l,l) + f(l+1,l+1) - f(l,l+1) )
$$
implies that \eqref{expr3} is equal to $-g(k_1,k_2)$.  \qed

\bigskip

{\it Proof of Lemma~\ref{odd-even}.} We use induction with respect to $n$. For $n=2$ the assertion is
easy to check. We assume that $n \ge 3$ and first consider the case that $n$ is odd. By the recursion,
it suffices to show that
\begin{multline}
\label{lr}
\sum_{ (l_1,\ldots, l_{(n-1)/2})}^{(k_1,\ldots,k_{(n-1)/2},k_{(n+1)/2})}
 \gamma(n-1,x;l_1,\ldots,l_{(n-1)/2}) \\ =
\sum_{ (l_1,\ldots, l_{(n-1)/2})}^{(k_1,\ldots,k_{(n-1)/2},2x+1-k_{(n+1)/2})}
 \gamma(n-1,x;l_1,\ldots,l_{(n-1)/2})
\end{multline}
By the induction hypothesis and by Lemma~\ref{hilf} (1) we know that
$$
\sum_{(l_{(n-1)/2})}^{(k''_{(n-1)/2},k_{(n+1)/2})} \gamma(n-1,x;l_1,\ldots,l_{(n-1)/2})
=
\sum_{(l_{(n-1)/2})}^{(k''_{(n-1)/2},2x+1-k_{(n+1)/2})} \gamma(n-1,x;l_1,\ldots,l_{(n-1)/2}).
$$
The assertion follows, since the left hand side of \eqref{lr} is equal to
$$
I_{k_{(n-1)/2}}^{k'_{(n-1)/2}, k''_{(n-1)/2}}
  \sum_{(l_1,\ldots,l_{(n-3)/2)})}^{(k_1,\ldots,k'_{(n-1)/2})}
 \sum_{(l_{(n-1)/2})}^{(k''_{(n-1)/2},k_{(n+1)/2})} \gamma(n-1,x;l_1,\ldots,l_{(n-1)/2})
$$
and the right hand side of \eqref{lr} is equal to
$$
I_{k_{(n-1)/2}}^{k'_{(n-1)/2}, k''_{(n-1)/2}}
  \sum_{(l_1,\ldots,l_{(n-3)/2)})}^{(k_1,\ldots,k'_{(n-1)/2})}
 \sum_{(l_{(n-1)/2})}^{(k''_{(n-1)/2},2x+1 - k_{(n+1)/2})} \gamma(n-1,x;l_1,\ldots,l_{(n-1)/2}).
$$
In the case that $n$ is even, Lemma~\ref{hilf} (2) is used in a similar way. \qed

\section{Derivation of the operator formula}
\label{derivation}

By Lemma~\ref{prop1}, we know that 
$
V_{k_i,k_{i+1}} \gamma(n,x;k_1,\ldots,k_{\lceil n/2 \rceil})
$
is antisymmetric in $k_i$ and $k_{i+1}$ for all $i$. Although the operators $V_{k_i,k_{i+1}}$
commute, this clearly does not imply that
$$
V_{k_1,k_2} V_{k_2,k_3} \dots V_{k_{\lceil n/2 \rceil -1},k_{\lceil n/2 \rceil}} 
\gamma(n,x;k_1,\ldots,k_{\lceil n/2 \rceil})
$$
is antisymmetric in $(k_1,k_2,\ldots,k_{\lceil n/2 \rceil})$. However, it is not hard to see that 
\begin{equation}
\label{star}
\left( \prod_{1 \le p < q \le \lceil n/2 \rceil} V_{k_p,k_q} \right) \gamma(n,x;k_1,\ldots,k_{\lceil n/2 \rceil})
\end{equation}
is antisymmetric in $(k_1,k_2,\ldots,k_{\lceil n/2 \rceil})$. This is a consequence of the following lemma, which generalizes \cite[Lemma~4]{fischer}. The proof is analogous to the
proof of  \cite[Lemma~4]{fischer} and thus we omit it here.

\begin{lem}
\label{product}
Let $W_{x,y}$ be an operator in $x$ and $y$, which is invertible as operator over 
$\mathbb{C}[x,y]$, and $W_{x_1,y_1} W_{x_2,y_2} = W_{x_2,y_2} W_{x_1,y_1}$
for all $x_1, x_2, y_1,y_2$.
Moreover, let $a(k_1,\ldots,k_{m})$ be a polynomial in $(k_1,\ldots,k_m)$. Then
$W_{k_i,k_{i+1}} a(k_1,\ldots,k_m)$ is antisymmetric in $k_i$ and $k_{i+1}$ for
all $i$ if and only if
$$
\left( \prod_{1 \le p < q \le m} W_{k_p,k_q} \right) a(k_1,\ldots,k_m)
$$
is antisymmetric in $(k_1,\ldots,k_m)$.
\end{lem}

We denote the polynomials in \eqref{star} by $\gamma^*(n,x;k_1,\ldots,k_{\lceil n/2 \rceil})$ and 
list them for $n=1,2,3,4,5$.
\begin{multline*}
1, x-k_1+1, \frac{1}{2} (2 x + 1 -k_1-k_2) (k_2-k_1), \frac{1}{6} (k_{2}-k_{1}) (2 x^3-3 k_{1} x^2-3 k_{2} x^2+6 x^2+k_{1}^2 x+k_{2}^2 x \\ -6 k_{1} x+4 k_{1} k_{2} x-6 k_{2}
   x+12 x+k_{1}^2-k_{1} k_{2}^2+k_{2}^2-6 k_{1}-k_{1}^2 k_{2}+4 k_{1} k_{2}-6 k_{2}+11), \\
\frac{1}{48} (k_{2}-k_{1}) (k_{3}-k_{1}) (k_{3}-k_{2}) (8 x^3-8 k_{1} x^2-8 k_{2} x^2-8 k_{3} x^2+12 x^2+2 k_{1}^2 x+2 k_{2}^2 x+2 k_{3}^2 x-8 k_{1} x \\ +6 k_{1} k_{2} x-8 k_{2} x+6 k_{1} k_{3} x+6 k_{2} k_{3} x-8 k_{3} x+30 x+k_{1}^2-k_{1}
k_{2}^2+k_{2}^2-k_{1} k_{3}^2  -k_{2} k_{3}^2+k_{3}^2-10 k_{1} \\ -k_{1}^2 k_{2}+3 k_{1} k_{2}-10 k_{2}-k_{1}^2 k_{3}-k_{2}^2 k_{3}+3
   k_{1} k_{3}-2 k_{1} k_{2} k_{3}+3 k_{2} k_{3}-10 k_{3}+25)
\end{multline*}
Although this list is shorter than the analog list for $\gamma(n,x;k_1,\ldots,k_{\lceil n/2 \rceil})$ 
(this is due to the factor $\prod\limits_{1 \le i < j \le \lceil n/2 \rceil} (k_j - k_i)$, which is
a consequence of the antisymmetry of the polynomial), it is still hard to guess the general pattern of 
$\gamma^*$. Thus we will apply a further operator to $\gamma^*$, in order to obtain a polynomial which 
factorizes into linear factors over $\mathbb{Q}$ and for which it is easy to recognize a pattern. 
This operator will have the property that it does not destroy the  
antisymmetry of the polynomial but restores the symmetry property of $\gamma$ given in Lemma~\ref{odd-even}. 
In the end, the fact that our operators are invertible will allow us to ``divide'' and give a formula 
for $\gamma$ itself.

The next lemma shows that the application of an operator, which is a symmetric polynomial in the 
shift operators, to an antisymmetric polynomial retains the antisymmetry.

\begin{lem}
\label{sym}
Let $a(k_1,\ldots,k_m)$ be a polynomial that is antisymmetric in $(k_1,\ldots,k_m)$ and $p(X_1,\ldots,X_m)$ be
a polynomial in $X_1,X_1^{-1},X_2,X_2^{-1},\ldots,X_m,X_m^{-1}$, which is symmetric in $(X_1,\ldots,X_m)$. 
Then $p(E_{k_1},\ldots,E_{k_m}) a(k_1,\ldots,k_m)$ is antisymmetric in $(k_1,\ldots,k_m)$.
\end{lem}

{\it Proof.} Let $\sigma \in {\mathcal S}_m$ be a permutation and
$p(X_1,\ldots,X_m) = \sum\limits_{(i_1,\ldots,i_m)} c_{i_1,\ldots,i_m} X_1^{i_1} \cdots X_m^{i_m}$.
The symmetry of $p(X_1,\ldots,X_m)$ implies that $c_{i_1,\ldots,i_m} = c_{i_{\sigma(1)},\ldots,i_{\sigma(m)}}$.
Thus
\begin{multline*}
p(E_{k_1},\ldots,E_{k_m}) a(k_1,\ldots,k_m)  =
\sum_{(i_1,\ldots,i_m)} c_{i_1,\ldots,i_m} E_{k_1}^{i_1}  \ldots E_{k_m}^{i_m} a(k_1,\ldots,k_m) \\ =
\sum_{(i_1,\ldots,i_m)} c_{i_1,\ldots,i_m} a(k_1+i_1,\ldots,k_m+i_m) \\
= \sgn \sigma \sum_{(i_1,\ldots,i_m)} c_{i_1,\ldots,i_m} a(k_{\sigma(1)}+i_{\sigma(1)},\ldots,k_{\sigma(m)}+i_{\sigma(m)}) \\
=  \sgn \sigma \sum_{(i_1,\ldots,i_m)} c_{i_{\sigma(1)},\ldots,i_{\sigma(m)}}
a(k_{\sigma(1)}+i_{\sigma(1)},\ldots,k_{\sigma(m)}+i_{\sigma(m)}) \\
=  \sgn \sigma \sum_{(i_1,\ldots,i_m)} c_{i_{1},\ldots,i_{m}}
a(k_{\sigma(1)}+i_{1},\ldots,k_{\sigma(m)}+i_{m}) \\
= \left. \sgn \sigma  \left( p(E_{l_1},\ldots,E_{l_m}) a(l_1,\ldots,l_m) \right) 
\right|_{(l_1,\ldots,l_m)=(k_{\sigma(1)},\ldots, k_{\sigma(m)})}.
 \qed
\end{multline*}

In the following lemma we identify operators whose application do not destroy 
symmetry properties of the type given in Lemma~\ref{odd-even}.

\begin{lem}
\label{negative}
Let $a(k_1,\ldots,k_m)$ be a polynomial such that 
$$a(k_1,\ldots,k_m)= \sigma \cdot a(k_1,\ldots,k_{m-1},d-k_m)$$ 
for $\sigma, d \in \mathbb{R}$ and 
$p(X_1,\ldots,X_m)$ be a polynomial in
$X_1,X_1^{-1},X_2,X_2^{-1},\ldots,X_m,X_m^{-1}$ such that 
$$p(X_1,\ldots,X_m)=p(X_1,\ldots,X_{m-1},X_m^{-1}).$$ 
Set $b(k_1,\ldots,k_m)=p(E_{k_1},\ldots,E_{k_m}) a(k_1,\ldots,k_m)$. Then we have 
$$b(k_1,\ldots,k_m)= \sigma \cdot b(k_1,\ldots,k_{m-1},d-k_m)$$
as well.
\end{lem}

{\it Proof.} Let $p(X_1,\ldots,X_m)=\sum\limits_{(i_1,\ldots,i_m)} c_{i_1,\ldots,i_m} 
X_1^{i_1} \cdots X_m^{i_m}$. By assumption $c_{i_1,\dots,i_{m-1},i_m}=c_{i_1,\cdots,i_{m-1},-i_m}$. Therefore,
\begin{multline*}
p(E_{k_1},\ldots,E_{k_m}) a(k_1,\ldots,k_m) = \sum_{(i_1,\ldots,i_m)} c_{i_1,\ldots,i_m} 
a(k_1+i_1,\ldots,k_m+i_m) \\
= \sigma \cdot \sum_{(i_1,\ldots,i_m)} c_{i_1,\ldots,i_m} a(k_1+i_1,\ldots,k_{m-1}+i_{m-1},d-k_m-i_m) \\
= \sigma \cdot \sum_{(i_1,\ldots,i_m)} c_{i_1,\ldots,i_{m-1},-i_m} a(k_1+i_1,\ldots,k_{m-1}+i_{m-1},d-k_m-i_m) \\
= \sigma \cdot \sum_{(i_1,\ldots,i_m)} c_{i_1,\ldots,i_{m}} a(k_1+i_1,\ldots,k_{m-1}+i_{m-1},d-k_m+i_m) \\
= \sigma \cdot \left. \left( p(E_{l_1},\ldots,E_{l_m}) a(l_1,\ldots,l_m) \right) \right|_{(l_1,\ldots,l_m)=(k_1,\ldots,k_{m-1},d-k_{m})}. \qed
\end{multline*}

\medskip

The previous two lemmas suggest to look for an operator $p(E_{k_1},\ldots,E_{k_{\lceil n/2 \rceil}})$, which 
is, on the one hand, symmetric in $(k_1,\ldots,k_{\lceil n/2 \rceil})$ and, on the other hand, has the 
property that the composition
$$
p(E_{k_1},\ldots,E_{k_{\lceil n/2 \rceil}}) \left( \prod_{1 \le p < q \le \lceil n/2 \rceil} (\id + E_{k_q} \Delta_{k_p}) \right)
$$
is invariant under the replacement of $E_{k_{\lceil n /2 \rceil}}$ by $E^{-1}_{k_{\lceil n /2 \rceil}}$.
This is accomplished in the following lemma.

\begin{lem}
\label{all}
The polynomial
$$
\left( \prod_{1 \le p < q \le \lceil n/2 \rceil} (\id + E_{k_q} \Delta_{k_p} ) E^{-1}_{k_p}
(\id  + E^{-1}_{k_q} \Delta_{k_p} ) \right)
\gamma(n,x;k_1,\ldots,k_{\lceil n/2 \rceil})
$$
is antisymmetric in $(k_1,\ldots,k_{\lceil n/2 \rceil})$. Moreover, if $n$ is odd then the polynomial is
invariant under the replacement of $k_i$ by $2x+1-k_i$ and if $n$ is even then the replacement of $k_i$ by
$2x+2-k_i$ only changes the sign of the polynomial.
\end{lem}

{\it Proof.} By Lemma~\ref{product},
$$
\left( \prod_{1 \le p < q \le \lceil n/2 \rceil } (\id + E_{k_q} \Delta_{k_p}) \right)
\gamma(n,x;k_1,\ldots,k_{\lceil n/2 \rceil})
$$
is antisymmetric in $(k_1,\ldots,k_{\lceil n/2 \rceil})$. Lemma~\ref{sym} and the fact that
$$
\prod_{1 \le p < q \le \lceil n/2 \rceil } X_p^{-1} (1 + X_q^{-1} (X_p - 1)) =
\prod_{1 \le p < q \le \lceil n/2 \rceil } (X_p^{-1} + X_q^{-1} - X_p^{-1} X_q^{-1} )
$$
is symmetric in $(X_1,\ldots,X_{\lceil n/2 \rceil})$ imply that the
expression in the statement of the lemma is still antisymmetric in $(k_1,\ldots,k_{\lceil n/2 \rceil})$. 

Next observe that the operator in the statement of the lemma is a polynomial in the shift 
operators $E_{k_i}^{\pm 1}$, which is invariant under the replacement of $E_{k_{\lceil n/2 \rceil}}$
by $E^{-1}_{k_{\lceil n/2 \rceil}}$.
Therefore, by
Lemma~\ref{odd-even} and by Lemma~\ref{negative}, the second assertion in lemma follows for $i=\lceil n/2 \rceil$.
The assertion for general $i$ follows from the antisymmetry of the polynomial. \qed

\medskip

The next lemma shows that the previous lemma together with the degree estimation
(Lemma~\ref{prop1}) determines $\gamma(n,x;k_1,\ldots,k_{\lceil n/2 \rceil})$ up to multiplicative
constant, which only depends on $n$.

\begin{lem}
\label{polynomial}
Let $p(k_1,\ldots,k_{\lceil n/2 \rceil})$ be an antisymmetric polynomial in
$(k_1,\ldots,k_{\lceil n/2 \rceil})$ over $\mathbb{C}$ of degree no greater than $n-1$ in every $k_i$
which is, in the case that $n$ is odd, invariant under the replacement of $k_i$ by $2x+1-k_i$ for every
$i$ and, in the case that $n$ is even, has the property that the replacement of $k_i$ by $2x+2-k_i$ only 
changes the sign of the
polynomial. Then $p(k_1,\ldots,k_{\lceil n/2 \rceil})$ equals
$$
C \cdot \prod_{1 \le i < j \le (n+1)/2} (k_j - k_i) (2x+1-k_i-k_j)
$$
if $n$ is odd and
$$
C \cdot \left( \prod_{1 \le i < j \le n/2} (k_j - k_i) (2x+2-k_i-k_j) \right) \prod_{i=1}^{n/2} (x+1-k_i)
$$
if $n$ is even, where $C$ is a constant in $\mathbb{C}$.
\end{lem}

{\it Proof.} We only consider the case that $n$ is even for the
other case is analogous. A polynomial $p(k_1,\ldots,k_{n/2})$ that
is antisymmetric in $(k_1,\ldots,k_{n/2})$ must have $k_j - k_i$
as a factor since
\begin{multline*}
p(k_1,\ldots,k_{i-1},k_i,k_{i+1},\ldots,k_{j-1},k_i,k_{j+1},\ldots, k_{n/2}) = \\
- p(k_1,\ldots,k_{i-1},k_i,k_{i+1},\ldots,k_{j-1},k_i,k_{j+1},\ldots, k_{n/2}).
\end{multline*}
This is because the polynomial changes the sign if we exchange the element in the $i$-th position
with the element in the $j$-th position.
If it furthermore has  the property that it will change the sign if $k_j$ is replaced by
$2x+2-k_j$ then the polynomial has a zero at $k_j=2x+2-k_i$ which explains the factor 
$2x+2-k_i-k_j$. Moreover it has a zero at
$k_i=x+1$ for every $i$, since
$$
p(k_1,\ldots,k_{i-1},x+1,k_{i+1},\ldots,k_{n/2}) = - p(k_1,\ldots,k_{i-1},x+1,k_{i+1},\ldots,k_{n/2}),
$$
which follows from $2x+2-(x+1)=x+1$. \qed

\medskip

Consequently, by Lemma~\ref{prop1}, Lemma~\ref{all} and Lemma~\ref{polynomial},
$$
\left( \prod_{1 \le p < q \le \lceil n/2 \rceil} (\id + E_{k_q} \Delta_{k_p} ) E^{-1}_{k_p}
(\id  + E^{-1}_{k_q} \Delta_{k_p} ) \right)
\gamma(n,x;k_1,\ldots,k_{\lceil n/2 \rceil})
$$
is equal to the polynomials given in Lemma~\ref{polynomial}. This determines $\gamma(n,x;k_1,\ldots,k_{\lceil n/2 \rceil})$
up to a multiplicative complex constant $C_n$. This is because the operators $\left(\id + E_{k_q} \Delta_{k_p}\right)$ and $\left(\id  + E^{-1}_{k_q} \Delta_{k_p}\right)$ are invertible by Lemma~\ref{invers}. In the following lemma we compute $C_n$.

\begin{lem}
\label{coeff}
 If $n$ is odd then 
$$C_n = \prod_{1 \le i < j \le (n+1)/2} \frac{1}{(j-i)(j+i-1)} = \frac{1}{(n-1)! (n-3)! \cdots 2!}$$ 
and if $n$ is even then
$$C_n = \prod_{1 \le i < j \le n/2} \frac{1}{(j-i)(j+i)} \prod_{i=1}^{n/2} \frac{1}{i} = \frac{1}{(n-1)! (n-3)! \cdots 1!}.$$
\end{lem}

{\it Proof.}  We expand the polynomial $\gamma(n,x;k_1,\ldots,k_{\lceil n/2 \rceil})$ with respect to the 
basis $\prod\limits_{i=1}^{\lceil n/2 \rceil} (k_i)_{m_i}$ and consider the coefficient of the 
basis element appearing in this expansion with maximal degree sequence  $(m_1,m_2,\ldots,m_{\lceil n/2 \rceil})$ in  lexicographic 
order. This coefficent is equal to $C_n$. We show by 
induction with respect to $n$ that this maximal degree sequence is $(n-1,n-3,\ldots,0)$ if 
$n$ is odd and $(n-1,n-3,\ldots,1)$ if $n$ is even. An analysis of the definition of 
$\sum\limits_{(l_1,\ldots,l_{m-1})}^{(k_1,\ldots,k_m)}$ and the induction hypothesis imply that  this maximal  
basis element of $\gamma(n,x;k_1,\ldots,k_{\lceil n/2 \rceil})$  is the maximal element of 
$$
\sum_{l_1=k_1}^{k_2-1} \sum_{l_2=k_2}^{k_3-1} \ldots \sum_{l_{(n-1)/2}=k_{(n-1)/2}}^{k_{(n+1)/2}}
\prod_{i=1}^{(n-1)/2} \frac{(l_i)_{n-2i}}{(n-2i)!}
$$
if $n$ is odd and the maximal element of 
$$
\sum_{l_1=k_1}^{k_2-1} \sum_{l_2=k_2}^{k_3-1} \ldots \sum_{l_{n/2}=k_{n/2}}^{x}
\prod_{i=1}^{n/2} \frac{(l_i)_{n-2i}}{(n-2i)!}
$$
if $n$ is even. \qed

\medskip

This immediately implies the following theorem.

\begin{theo}
\label{sidemain}
If $n$ is odd then
\begin{multline*}
\gamma(n,x;k_1,\ldots,k_{(n+1)/2}) \\
= \left( \prod_{1 \le p < q \le (n+1)/2} (\id + E_{k_q} \Delta_{k_p})^{-1} E_{k_p} (\id + E_{k_q}^{-1} \Delta_{k_p})^{-1} \right) \\
\prod_{1 \le i < j \le (n+1)/2} \frac{(k_j - k_i) (2x+1-k_i-k_j)}{(j-i)(j+i-1)}
\end{multline*}
and if $n$ is even then
\begin{multline*}
\gamma(n,x;k_1,\ldots,k_{n/2}) \\
= \left( \prod_{1 \le p < q \le n/2} (\id + E_{k_q} \Delta_{k_p})^{-1} E_{k_p} (\id + E_{k_q}^{-1} \Delta_{k_p})^{-1} \right) \\
\prod_{1 \le i < j \le n/2} \frac{(k_j-k_i) (2x+2-k_i-k_j)}{(j-i)(j+i)}  \prod_{i=1}^{n/2} \frac{x+1-k_i}{i}.
\end{multline*}
\end{theo}

{\it Proof.} The assertion follows from Lemma~\ref{invers}, Lemma~\ref{prop1}, Lemma~\ref{all}, Lemma~\ref{polynomial} and Lemma~\ref{coeff}.  \qed

\bigskip

We are finally able to prove the main theorem.

\medskip

{\it Proof of Theorem~\ref{main}.} Observe that
$$
\prod_{1 \le p <  q \le \lceil n/2 \rceil} (1+X_q (X_p -1))(1+X_q^{-1} (X_p-1))(1+X_q (X_p^{-1}-1)) (1+X_q^{-1}(X_p^{-1}-1))
$$
is invariant under the replacement of $X_i$ by $X^{-1}_i$. Moreover, it is symmetric in $(X_1,\ldots,X_{\lceil n/2})$, since
the factor associated to the pair $(p,q)$ is equal to
\begin{multline*}
(1+X_p X_q - X_q) (1+X_p^{-1} X_q^{-1} - X_q^{-1}) X_p X_q \\
\times (1+X_p X_q^{-1} - X_q^{-1}) X_q \\
\times (1+X_p^{-1} X_q - X_q) X_p \\
\times X_p^{-2} X_q^{-2}
\end{multline*}
and this is symmetric in $X_p$ and $X_q$ as every line is.
Thus, by Lemma~\ref{sym},
\begin{multline}
\label{odd1}
\left( \prod_{1 \le p < q \le (n+1)/2} (\id+E_{k_q} \Delta_{k_p})(\id+E_{k_q}^{-1} \Delta_{k_p})
(\id - E_{k_q} E_{k_p}^{-1} \Delta_{k_p}) (\id- E^{-1}_{k_q} E_{k_p}^{-1} \Delta_{k_p}) \right) \\
\prod_{1 \le i < j \le (n+1)/2} (k_j - k_i) (2x+1-k_i-k_j)
\end{multline}
is antisymmetric in $(k_1,\ldots,k_{(n+1)/2})$ if $n$ is odd and
\begin{multline}
\label{even1}
\left( \prod_{1 \le p <  q \le n/2} (\id+E_{k_q} \Delta_{k_p})(\id+E_{k_q}^{-1} \Delta_{k_p})
(\id - E_{k_q} E_{k_p}^{-1} \Delta_{k_p}) (\id - E^{-1}_{k_q} E_{k_p}^{-1} \Delta_{k_p}) \right) \\
\prod_{1 \le i < j \le n/2} (k_j-k_i) (2x+2-k_i-k_j)  \prod_{i=1}^{n/2} (x+1-k_i)
\end{multline}
is antisymmetric in $(k_1,\ldots,k_{n/2})$ if $n$ is even. Moreover, by Lemma~\ref{negative},
\eqref{odd1} is invariant
under the replacement of $k_i$ by $2x+1-k_i$, whereas \eqref{even1} changes the sign if $k_i$ is
replaced by $2x+2-k_i$. Consequently, by Lemma~\ref{polynomial}, \eqref{odd1} is equal to
$$
D_n \cdot \prod_{1 \le i < j \le (n+1)/2} (k_j - k_i) (2x+1-k_i-k_j)
$$
if $n$ is odd and \eqref{even1} is equal to
$$
D_n \cdot \prod_{1 \le i < j \le n/2} (k_j - k_i) (2x+2-k_i-k_j)  \prod_{i=1}^{n/2} (x+1-k_i)
$$
if $n$ is even.  If we compare the coefficient of a monomial of
maximal degree we see that $D_n=1$. Now, if $n$ is odd
then, by Theorem~\ref{sidemain},
\begin{multline*}
\gamma(n,x;k_1,\ldots,k_{(n+1)/2}) \\
= \left( \prod_{1 \le p < q \le (n+1)/2} (\id + E_{k_q} \Delta_{k_p})^{-1} E_{k_p} (\id + E_{k_q}^{-1} \Delta_{k_p})^{-1} \right) \\
\prod_{1 \le i < j \le (n+1)/2} \frac{(k_j - k_i) (2x+1-k_i-k_j)}{(j-i)(j+i-1)}  \\
= \left( \prod_{1 \le p < q \le (n+1)/2} (\id + E_{k_q} \Delta_{k_p})^{-1} E_{k_p} (\id + E_{k_q}^{-1} \Delta_{k_p})^{-1} \right) \\
 \left( \prod_{1 \le p <  q \le (n+1)/2} (\id+E_{k_q} \Delta_{k_p})(\id+E_{k_q}^{-1} \Delta_{k_p})
(\id - E_{k_q} E_{k_p}^{-1} \Delta_{k_p}) (\id- E^{-1}_{k_q} E_{k_p}^{-1} \Delta_{k_p})  \right) \\
 \prod_{1 \le i < j \le (n+1)/2} \frac{(k_j - k_i) (2x+1-k_i-k_j)}{(j-i)(j+i-1)}  \\
= 
\left( \prod_{1 \le p < q \le (n+1)/2}
E_{k_p} (\id - E_{k_q} E_{k_p}^{-1} \Delta_{k_p}) (\id- E^{-1}_{k_q} E_{k_p}^{-1} \Delta_{k_p}) \right) \\
\prod_{1 \le i < j \le (n+1)/2} \frac{(k_j - k_i) (2x+1-k_i-k_j)}{(j-i)(j+i-1)}
\end{multline*}
Similarly, if $n$ is even then
\begin{multline*}
\gamma(n,x;k_1,\ldots,k_{n/2}) \\
= 
\left( \prod_{1 \le p < q \le n/2}
E_{k_p} (\id - E_{k_q} E_{k_p}^{-1} \Delta_{k_p}) (\id- E^{-1}_{k_q} E_{k_p}^{-1} \Delta_{k_p}) \right) \\
\prod_{1 \le i < j \le n/2} \frac{(k_j - k_i) (2x+2-k_i-k_j)}{(j-i)(j+i)}
\prod_{i=1}^{n/2} \frac{x+1-k_i}{i}. \qed
\end{multline*}

\section{From operator formulas to generating functions}
\label{generating}

In this section we follow a hint of Doron Zeilberger and translate the operator formulas into generating function results. We start out with ordinary
monotone triangles. In~\cite{fischer} we have shown that the number of monotone
triangles with bottom row $(k_1,\ldots,k_n)$ is given by
$$
\alpha(n;k_1,\ldots,k_n) = \left( \prod_{1 \le p < q \le n} (\id +
E_{k_p} E_{k_q} - E_{k_p}) \right) \prod_{1 \le i < j \le n}
\frac{k_j - k_i}{j-i}
$$
if $k_1 < k_2 < \dots < k_n$ and $k_i \in \mathbb{Z}$.
We define $\overline{\alpha}_c(n;k_1,\ldots,k_n)=\prod\limits_{1
\le i < j \le n} \frac{k_j - k_i}{j-i}$ if $k_l \ge c$ for all $l$
and zero elsewhere. Then
$$
\alpha(n;k_1,\ldots,k_n) = \left( \prod_{1 \le p < q \le n} (\id +
E_{k_p} E_{k_q} - E_{k_p}) \right)
\overline{\alpha}_c(n;k_1,\ldots,k_n)
$$
for all $(k_1,\ldots,k_n)$ with $k_l \ge c$. We compute the
generating function
\begin{equation}
\label{genfun} \sum_{(k_1,\ldots,k_n) \ge (-n+1,\ldots,-n+1)}
X_1^{k_1}  \ldots X_n^{k_n} \left( \prod_{1 \le p < q \le n} (\id
+ E_{k_p} E_{k_q} - E_{k_p}) \right)
\overline{\alpha}_0(n;k_1,\ldots,k_n).
\end{equation}
Thus, the coefficient of $X_1^{k_1} X_2^{k_2} \dots X_n^{k_n}$
gives the number of monotone triangles with bottom row
$(k_1,\ldots,k_n)$ if $0 \le k_1 < k_2 < \ldots < k_n$ and $k_i \in \mathbb{Z}$. 
Let
$$
\prod_{1 \le p < q \le n} (1 + Y_p Y_q - Y_p) =
\sum_{(j_1,\ldots,j_n) \atop 0 \le j_i \le n-1} a(j_1,\ldots,j_n)
Y_1^{j_1} \cdots Y_n^{j_n}.
$$
Using this notation, \eqref{genfun} is equal to
\begin{equation}
\label{genfun1} \sum_{(k_1,\ldots,k_n) \ge (-n+1,\ldots,-n+1)}
\sum_{(j_1,\ldots,j_n) \atop 0 \le j_i \le n-1} \,
a(j_1,\ldots,j_n) X_1^{k_1}  \ldots X_n^{k_n} \,
\overline{\alpha}_0(n;k_1+j_1,\ldots,k_n+j_n)
\end{equation}
We set $(l_1,\ldots,l_n)=(k_1+j_1,\ldots,k_n+j_n)$. Consequently,
\eqref{genfun1} is equal to
\begin{equation}
\label{genfun2} \sum_{(j_1,\ldots,j_n) \atop 0 \le j_i \le n-1}
a(j_1,\ldots,j_n) X_1^{-j_1}  \ldots X_n^{-j_n}
\sum_{(l_1,\ldots,l_n) \ge (-n+1+j_1,\ldots,-n+1+j_n)}
\overline{\alpha}_0(n;l_1,\ldots,l_n) X_1^{l_1} \dots X_n^{l_n}.
\end{equation}
Since $\overline{\alpha}_0(l_1,\ldots,l_n)=0$ if $l_i <0$ for an
$i$ and $j_l \le n-1$ for all $l$, \eqref{genfun2} is equal to
\begin{multline*}
\sum_{(j_1,\ldots,j_n) \atop 0 \le j_i \le n-1} a(j_1,\ldots,j_n)
X_1^{-j_1}  \ldots X_n^{-j_n} \sum_{(l_1,\ldots,l_n) \ge
(0,\ldots,0)} \overline{\alpha}_0(n;l_1,\ldots,l_n) X_1^{l_1}
\dots X_n^{l_n} \\
= \prod_{1 \le p < q \le n} (1+X_p^{-1} X_q^{-1} - X_p^{-1})
\sum_{(l_1,\ldots,l_n) \ge (0,\ldots,0)} X_1^{l_1} \dots X_n^{l_n}
\prod_{1 \le i < j \le n} \frac{l_j - l_i}{j-i}.
\end{multline*}
The Vandermonde determinant evaluation implies that 
$$
\prod_{1 \le i < j \le n} \frac{l_j - l_i}{j-i} = \det_{1 \le i, j
\le n} \binom{l_i}{j-1}
$$
and, consequently, the generating function is equal to
$$
\prod_{1 \le p < q \le n} (1+X_p^{-1} X_q^{-1} - X_p^{-1})
\sum_{(l_1,\ldots,l_n) \ge (0,\ldots,0)} \det_{1 \le i, j \le n}
X_i^{l_i} \binom{l_i}{j-1}.
$$
Observe that
\begin{multline*}
\sum_{l=0}^{\infty} X^l \binom{l}{j-1} = \frac{X^{j-1}}{(j-1)!}
\sum_{l=0}^{\infty} l (l-1) \dots (l-j+2) X^{l-j+1} =
\frac{X^{j-1}}{(j-1)!} \frac{d}{d X^{j-1}} \left(
\sum_{l=0}^{\infty} X^l \right) \\ = \frac{X^{j-1}}{(j-1)!}
\frac{d}{d X^{j-1}} (1-X)^{-1} = \frac{X^{j-1}}{(1-X)^j}.
\end{multline*}
Therefore, the generating function is equal to
\begin{multline*}
\prod_{1 \le p < q \le n} (1+X_p^{-1} X_q^{-1} - X_p^{-1})
\frac{1}{(1-X_1)(1-X_2)\dots(1-X_n)} \det_{1 \le i, j \le n}
\left( \frac{X_i}{1-X_i} \right)^{j-1} \\
 = \frac{1}{(1-X_1)(1-X_2)\dots(1-X_n)} \prod_{1 \le i < j \le n}
(1+X_i^{-1} X_j^{-1} - X_i^{-1}) \left(
\frac{X_j}{(1-X_j)}-\frac{X_i}{(1-X_j)} \right) \\
= \frac{1}{(1-X_1)(1-X_2)\dots(1-X_n)} \prod_{1 \le i < j \le n}
\frac{(X_j - X_i)(1-X_j+X_i X_j)}{X_i (1-X_i) X_j (1-X_j)} \\ =
\prod_{i=1}^n \frac{1}{X_i^{n-1} (1-X_i)^n} \prod_{1 \le i < j \le
n} (X_j - X_i)(1-X_j+X_i X_j),
\end{multline*}
where the Vandermonde determinant evaluation is used again. We summerize the result in the following theorem.

\begin{theo}
\label{theogf1}
The coefficient of $X_1^{k_1} X_2^{k_2} \dots X_n^{k_n}$ in 
\begin{equation}
\label{gf1}
\prod_{i=1}^n \frac{1}{X_i^{n-1} (1-X_i)^n} \prod_{1 \le i < j \le
n} (X_j - X_i)(1-X_j+X_i X_j)
\end{equation}
is equal to the number of monotone triangles with bottom row $(k_1,\ldots,k_n)$ if 
$0 \le k_1 < k_2 < \ldots < k_n$ and $k_i \in \mathbb{Z}$, where \eqref{gf1} is interpreted 
as a formal laurent series, $\frac{1}{1-X_i} = \sum\limits_{j=0}^{\infty} X_i^j$.
\end{theo}

Consequently, the enumeration of $n \times n$ alternating sign matrices amounts to 
compute the constant term of 
$$
\prod_{i=1}^n \frac{1}{X_i^{n+i-2} (1-X_i)^n} \prod_{1 \le i < j \le
n} (X_j - X_i)(1-X_j+X_i X_j).
$$
This is because monotone triangles with bottom row $(0,1,\ldots,n-1)$ correspond to 
$n \times n$ alternating sign matrices. Zeilberger~\cite{zeilberger} has 
used constant term identities to give the first proof of the alternating 
sign matrices theorem. (His identities are different from our result.)

Note that, if, for instant, we choose $n=3$ in the generating function 
in Theorem~\ref{theogf1} and consider the coefficient of $X_1^3 X_2^2 X_3$ we
obtain $-1$, which is oviously not the number of monotone triangles 
with bottom row $(3,2,1)$, since there exists no monotone triangle with this 
property. This coefficient is of course the values of $\alpha(3;3,2,1)$. On the other hand, if we 
consider monomials $X_1^{k_1} \cdots X_n^{k_n}$ with negative exponents then their
coefficients are not equal to $\alpha(n;k_1,\ldots,k_n)$: for example the coefficient of $X_1^{-1} X_2^2 X_3^3$ 
is $7$ and this is  not  $\alpha(3;-1,2,3)=23$. (In order to compute the number of monotone 
triangles with bottom row $(-1,2,3)$ using the generating function from Theorem~\ref{theogf1}, one 
can make use of the fact that 
$\alpha(3;-1,2,3)=\alpha(3;-1+c,2+c,3+c)$ for all integers $c$.)

\bigskip

Next we derive an analog generating function for halved monotone
triangles with prescribed bottom row. Observe that
Theorem~\ref{main} is equivalent to
\begin{multline*}
\gamma(n,x;k_1,\ldots,k_{\lceil n/2 \rceil}) =  \\ \left( \prod_{1
\le p < q \le \lceil n/2 \rceil} (E^{-1}_{k_p} + E^{-1}_{k_q} -
\id) (E^{-1}_{k_p} E^{-1}_{k_q} + \id - E^{-1}_{k_q} ) \right)
\overline{\gamma}(n,x;k_1,\ldots,k_{\lceil n/2 \rceil})
\end{multline*}
where
$$
\overline{\gamma}(n,x;k_1,\ldots,k_{(n+1)/2}) =
\prod_{1 \le i < j \le
(n+1)/2} \frac{(k_j - k_i) (2x+2-n-k_i-k_j)}{(j-i)(j+i-1)}
$$
if $n$ is odd and
$$
\overline{\gamma}(n,x;k_1,\ldots,k_{n/2}) = 
 \prod_{1 \le i < j \le n/2} \frac{(k_j - k_i)
(2x+4-n-k_i-k_j)}{(j-i)(j+i)} \prod_{i=1}^{n/2} \frac{x+2-n/2-k_i}{i}
$$
if $n$ is even. Here, we define
$\overline{\gamma}_c(n,x;k_1,\ldots,k_{\lceil n/2 \rceil})$ to be
equal to $\overline{\gamma}(n,x;k_1,\ldots,k_{\lceil n/2 \rceil})$
if $k_l \le c$ for all $l$ and zero elsewhere. Then
\begin{multline*}
\gamma(n,x;k_1,\ldots,k_{\lceil n/2 \rceil}) = \\ \left( \prod_{1
\le p < q \le \lceil n/2 \rceil} (E^{-1}_{k_p} + E^{-1}_{k_q} -
\id) (E^{-1}_{k_p} E^{-1}_{k_q} + \id - E^{-1}_{k_q} ) \right)
\overline{\gamma}_c(n,x;k_1,\ldots,k_{\lceil n/2 \rceil})
\end{multline*}
for all $(k_1,\ldots,k_{\lceil n/2 \rceil})$ with $k_l \le c$. We
compute the generating function
\begin{multline}
\label{genfunhalved} \sum_{(k_1,\ldots,k_{\lceil n/2 \rceil}) \le
(n-1+c,\ldots,n-1+c)} X_1^{k_1}  \ldots X_{\lceil n/2
\rceil}^{k_{\lceil n/2 \rceil}}  \\ \times \left( \prod_{1 \le p <
q \le \lceil n/2 \rceil } (E^{-1}_{k_p} + E^{-1}_{k_q} - \id) (E^{-1}_{k_p} E^{-1}_{k_q} +\id -
E^{-1}_{k_q})  \right)
\overline{\gamma}_c(n,x;k_1,\ldots,k_{\lceil n/2 \rceil}).
\end{multline}
Let
$$
\prod_{1 \le p < q \le \lceil n/2 \rceil} (Y_p^{-1} + Y_q^{-1} -
1) (Y_p^{-1} Y_q^{-1} + 1 -Y_q^{-1}) =
\sum_{(j_1,\ldots,j_{\lceil n/2 \rceil}) \atop -n+1 \le j_i \le 0}
b(j_1,\ldots,j_{\lceil n/2 \rceil}) Y_1^{j_1} \cdots Y_{\lceil n/2
\rceil}^{j_{\lceil n/2 \rceil}}.
$$
Thus, the generating function \eqref{genfunhalved} is equal to
\begin{multline}
\label{genfunhalved1} \sum_{(k_1,\ldots,k_{\lceil n/2 \rceil}) \le
(n-1+c,\ldots,n-1+c)}  \sum_{(j_1,\ldots,j_{\lceil n/2 \rceil})
\atop -n+1 \le j_i \le 0}  b(j_1,\ldots,j_{\lceil n/2 \rceil})
X_1^{k_1}  \ldots X_{\lceil n/2 \rceil}^{k_{\lceil n/2 \rceil}} \\
\times \overline{\gamma}_c(n,x;k_1+j_1,\ldots,k_{\lceil n/2
\rceil}+j_{\lceil n/2 \rceil}).
\end{multline}
Again we set $(l_1,\ldots,l_{\lceil n/2
\rceil})=(k_1+j_1,\ldots,k_{\lceil n/2 \rceil}+j_{\lceil n/2
\rceil})$. Consequently, \eqref{genfunhalved1} is equal to
\begin{multline}
\label{genfunhalved2} \sum_{(j_1,\ldots,j_{\lceil n/2 \rceil})
\atop -n+1 \le j_i \le 0}  b(j_1,\ldots,j_{\lceil n/2 \rceil})
X_1^{-j_1}  \ldots X_{\lceil n/2 \rceil}^{-j_{\lceil n/2 \rceil}}
\\ \times \sum_{(l_1,\ldots,l_{\lceil n/2 \rceil}) \le (n-1+j_1+c,\ldots,n-1+j_{\lceil n/2 \rceil}+c)}
\overline{\gamma}_c(n,x;l_1,\ldots,l_{\lceil n/2 \rceil})
X_1^{l_1} \ldots X_{\lceil n/2 \rceil}^{l_{\lceil n/2 \rceil}}.
\end{multline}
Since $j_l \ge -n+1$ for all $l$ and
$\overline{\gamma}_c(n,x;l_1,\ldots,l_{\lceil n/2 \rceil})=0$ if
$l_i > c$ for an $i$, \eqref{genfunhalved2} is equal to
\begin{multline}
\label{genfunhalved3} \sum_{(j_1,\ldots,j_{\lceil n/2 \rceil})
\atop -n+1 \le j_i \le 0}  b(j_1,\ldots,j_{\lceil n/2 \rceil})
X_1^{-j_1}  \ldots X_{\lceil n/2 \rceil}^{-j_{\lceil n/2 \rceil}}
\\ \times \sum_{(l_1,\ldots,l_{\lceil n/2 \rceil}) \le (c,\ldots,c)}
\overline{\gamma}_c(n,x;l_1,\ldots,l_{\lceil n/2 \rceil})
X_1^{l_1}  \ldots X_{\lceil n/2 \rceil}^{l_{\lceil n/2 \rceil}} \\
= \prod_{1 \le p < q \le \lceil n/2 \rceil} (X_p + X_q -
1) (X_p X_q + 1 - X_q) \\
\times \sum_{(l_1,\ldots,l_{\lceil n/2 \rceil}) \le (c,\ldots,c)}
\overline{\gamma}_c(n,x;l_1,\ldots,l_{\lceil n/2 \rceil})
X_1^{l_1} \ldots X_{\lceil n/2 \rceil}^{l_{\lceil n/2 \rceil}}.
\end{multline}
In the last expression $\overline{\gamma}_c$ can be replaced by
$\overline{\gamma}$. The following lemma provides us with 
determinantal expressions for
$\overline{\gamma}(n,x;k_1,\ldots,k_{\lceil n/2 \rceil})$.

\begin{lem}
\label{determinants} (1) $$\det_{1 \le i, j \le n}
\binom{k_i+j-1}{2j-1} = 
\prod_{1 \le i < j \le n} \frac{(k_j - k_i)(k_i+k_j)}{(j-i)(j+i)} \prod_{i=1}^{n}
\frac{k_i}{i}$$

(2) $$\det_{1 \le i, j \le n} \binom{k_i+j-3/2}{2j-2} =
\prod_{1 \le i < j \le n} \frac{(k_j -
k_i)(k_i+k_j)}{(j-i)(j+i-1)}$$
\end{lem}

 {\it Proof.} We only prove (1) since the proof of (2) is
 similar. First observe that
 $$
 \binom{k_i+j-1}{2j-1} = \frac{k_i}{(2j-1)!} \prod_{l=1}^{j-1} (k_i^2-l^2).
 $$
 Thus, the determinant in (1) is equal to
 $$
\prod_{i=1}^{n-1} \frac{1}{(2i+1)!} \prod_{i=1}^{n} k_i \det_{1
\le i, j \le n} \left( \prod_{l=1}^{j-1} (k_i^2 - l^2) \right).
$$
The assertion follows from
$$
\det_{1 \le i, j \le n}  p_j(Y_i) = \prod_{1 \le i < j \le n} (Y_j - Y_i),
$$
where $p_j(Y)$ is a polynomial in $Y$ of degree $j-1$ whose leading coefficient is $1$.
This is a consequence of the Vandermonde determinant evaluation.
 \qed

\medskip

Lemma~\ref{determinants} implies that
\begin{equation}
\label{det1}
 \overline{\gamma}(n,x;k_1,\ldots,k_{(n+1)/2}) =
 (-1)^{\binom{(n+1)/2}{2}} \det_{1 \le i,j \le (n+1)/2}
 \binom{k_i+j+n/2-x-5/2}{2j-2}
\end{equation}
 if $n$ is odd and
\begin{equation}
\label{det2}
\overline{\gamma}(n,x;k_1,\ldots,k_{n/2}) =
 (-1)^{\binom{(n+2)/2}{2}} \det_{1 \le i,j \le n/2}
 \binom{k_i+j+n/2-x-3}{2j-1}
\end{equation}
if $n$ is even.
If we use these determinantal presentations for $\overline{\gamma}$ in \eqref{genfunhalved3}, we 
obtain the following generating function 
\begin{multline}
\label{odd}
\prod_{1 \le p < q \le (n+1)/2} (1-X_p - X_q)(X_p X_q+1-X_q) \\
\times \det_{1 \le i,j \le (n+1)/2} 
\left( \sum_{l_i=-\infty}^c \binom{l_i+j+n/2-x-5/2}{2j-2} X_i^{l_i} \right)
\end{multline}
if $n$ is odd, and 
\begin{multline}
\label{even}
\prod_{1 \le p < q \le n/2} (1-X_p - X_q)(X_p X_q+1-X_q) \\
\times \det_{1 \le i,j \le n/2} 
\left(- \sum_{l_i=-\infty}^c \binom{l_i+j+n/2-x-3}{2j-1} X_i^{l_i} \right)
\end{multline}
if $n$ is even.

If we choose $c=x+1/2-n/2$ in case that $n$ is odd and $c=x+2-n/2$ in case that $n$ is even, 
the determinants in the expression 
above simplify. This follows from the following identities.  
\begin{multline}
\label{oddid}
\sum_{l=-\infty}^c \binom{l+j-c+z}{2j-2} X^l =
\frac{X^{j+c-z-2}}{(2j-2)!} \frac{d}{d X^{2j-2}} \left(
\sum_{l=-\infty}^c
X^{l+j-c+z} \right) \\
= \frac{X^{j+c-z-2}}{(2j-2)!} \frac{d}{d X^{2j-2}} \left(
\frac{X^{j+z+1}}{X-1} \right) =  \frac{X^{j+c-z-2}}{(2j-2)!} \frac{d}{d
X^{2j-2}} \left( \frac{1}{X-1} \right) =
\frac{X^{j+c-z-2}}{(X-1)^{2j-1}}
\end{multline}
where $z=-j-1,-j,\ldots,j-3$. (Note that the identity is true for all $j$ if $z=-2$.)
Similarly,
\begin{equation}
\label{evenid}
\sum\limits_{l=-\infty}^c \binom{l+j-c+z}{2j-1} X^l = -
\frac{X^{j+c-z-1}}{(X-1)^{2j}}
\end{equation}
where $z=-j-1,-j,\ldots,j-2$. (This identity is true for all $j$ if $z=-2,-1$.)

We first consider the case that $n$ is
odd. By \eqref{oddid} ($z=-2$) the generating function \eqref{odd} is equal to
\begin{multline*}
\prod_{1 \le p < q \le (n+1)/2 } (1- X_p - X_q) (X_p X_q + 1 - X_q)
 \det_{1 \le i,j \le (n+1)/2} \left(
\frac{X_i^{j+c}}{(X_i-1)^{2j-1}} \right) \\
= \prod_{1 \le p < q \le (n+1)/2 } (1- X_p - X_q) (X_p X_q + 1 - X_q) \prod_{i=1}^{(n+1)/2} \frac{X_i^{c+1}}{(X_i-1)} \det_{1 \le
i,j \le (n+1)/2} \left( \frac{X_i}{(X_i-1)^2} \right)^{j-1}. 
\end{multline*}
The Vandermonde determinant evaluation shows that this is equal to 
\begin{multline}
 \label{endeodd}
 \prod_{1 \le p < q \le (n+1)/2 } (1- X_p - X_q) (X_p X_q + 1 - X_q) \prod_{i=1}^{(n+1)/2} \frac{X_i^{c+1}}{(X_i-1)} \\ \times
 \prod_{1 \le i < j \le (n+1)/2} \left( \frac{X_j}{(X_j-1)^2} -
\frac{X_i}{(X_i-1)^2} \right) \\
= \prod_{1 \le i < j \le (n+1)/2} (X_j-X_i)(X_i+X_j-1)(X_i X_j -
1)(1- X_j + X_i X_j ) \prod_{i=1}^{(n+1)/2}
\frac{X_i^{x+3/2-n/2}}{(X_i-1)^n}.
\end{multline}

Finally we consider the case that $n$ is even. 
By \eqref{evenid} ($z=-1$) the generating function \eqref{even} is equal to
\begin{multline*}
  \prod_{1 \le p < q \le n/2 } (1- X_p - X_q) (X_p X_q + 1 - X_q)
 \det_{1 \le i,j \le n/2} \left( \frac{X_i^{j+c}}{(X_i-1)^{2j}}
 \right) = \\
 \prod_{1 \le p < q \le n/2 } (1- X_p - X_q) (X_p X_q + 1 - X_q)
\prod_{i=1}^{n/2} \frac{X_i^{c+1}}{(X_i-1)^2}
 \det_{1 \le i,j \le n/2} \left( \frac{X_i}{(X_i-1)^{2}}
 \right)^{j-1} 
\end{multline*}
The Vandermonde determinant evaluation now shows that this is equal to
\begin{multline}
\label{endeeven}
\prod_{1 \le p < q \le n/2 } (1- X_p - X_q) (X_p X_q + 1 - X_q)
\prod_{i=1}^{n/2} \frac{X_i^{c+1}}{(X_i-1)^2} \\
 \times \prod_{1 \le i < j \le n/2}  \left( \frac{X_j}{(X_j-1)^2} -
\frac{X_i}{(X_i-1)^2} \right) \\
= \prod_{1 \le i < j \le n/2} (X_j-X_i)(X_i+X_j-1)(X_i X_j -
1)(1-X_j+X_i X_j) \prod_{i=1}^{n/2} \frac{X_i^{x+3-n/2}}{(X_i-1)^n}.
\end{multline}

In this case the generating functions in \eqref{endeodd} and \eqref{endeeven}
are understood as formal laurent series in $1/X_i$, i.e. 
$1/(X_i-1) = \sum_{j=-\infty}^{-1} X_i^j$.

\end{document}